\documentclass[a4paper, 12pt]{amsart}

\usepackage{amsmath, amsthm}

\theoremstyle{plain}
\newtheorem{thm}{Theorem}[section]
\newtheorem{lemma}[thm]{Lemma}
\newtheorem{cor}[thm]{Corollary}
\newtheorem{conjecture}[thm]{Conjecture}

\theoremstyle{remark}
\newtheorem{example}[thm]{Example}

\providecommand{\aut}{\mathop{\rm Aut \,}\nolimits}

\providecommand{\sym}{\mathop{\rm Sym \,}\nolimits}
\providecommand{\Wr}{\mathop{\rm Wr \,}\nolimits}

\renewcommand{\\}{\vspace{3mm}}

\begin{document}

\title[\bf Subdegree growth rates of primitive groups]{\bf Subdegree growth rates of infinite primitive permutation groups}

\author{{\bf Simon M. Smith}}

\address{Mathematical Institute, University of Oxford}

\email{simon.smith@chch.oxon.org}

\subjclass{20 B 15, 05 C 25}

\date{\today}

\begin{abstract}

A transitive group $G$ of permutations of a set $\Omega$ is primitive if the only  
$G$-invariant equivalence relations on $\Omega$ are the trivial and universal relations.

If $\alpha \in \Omega$, then the orbits of the stabiliser $G_\alpha$ on $\Omega$ are called  
the $\alpha$-suborbits of $G$; when $G$ acts transitively the cardinalities of these
$\alpha$-suborbits are the subdegrees of $G$.

If $G$ acts primitively on an infinite set $\Omega$, and all the suborbits of $G$
are finite, Adeleke and Neumann asked if, after enumerating the subdegrees of $G$ as
a non-decreasing sequence $1 = m_0 \leq m_1 \leq \cdots$, the subdegree growth
rates of infinite
primitive groups that act distance-transitively on locally finite
distance-transitive graphs are extremal, and conjecture there might
exist a number $c$ which perhaps depends upon $G$, perhaps only on
$m$, such that $m_r \leq c(m-2)^{r-1}$.

In this paper it is shown that such an enumeration is not desirable, as
there exist infinite primitive permutation groups possessing no
infinite subdegree, in which two distinct subdegrees are each equal to the cardinality of  
infinitely many suborbits. The examples used to show this provide several novel methods for constructing
infinite primitive graphs.

A revised enumeration method is then proposed, and it is shown that, under
this, Adeleke and Neumann's question may be answered, at least for groups  exhibiting suitable rates of growth.

\end{abstract}

\maketitle

\section{Introduction}

Let $G$ be a group of permutations of an
infinite set $\Omega$. If $\alpha \in \Omega$ and $g \in G$, we denote
the image of $\alpha$ under $g$ by $\alpha^g$, and the set of images
of $\alpha$ under all elements of $G$ by $\alpha^G$. All permutations
will act on the right. We denote the stabiliser of $\alpha$ in $G$ by
$G_\alpha$, and if $\Sigma \subseteq \Omega$ we denote the setwise and
pointwise stabilisers of $\Sigma$ in $G$ by $G_{\{\Sigma\}}$ and
$G_{(\Sigma)}$ respectively. A transitive group $G$ is {\it primitive}
on $\Omega$ if the only $G$-invariant equivalence relations on
$\Omega$ are the trivial and universal relations. It is well known
that all point stabilisers in a primitive group are maximal.

Given $\alpha, \beta \in \Omega$, the set $(\alpha, \beta)^G$ is
called an {\it orbital} of $G$. It is {\it diagonal} if $\alpha$ and
$\beta$ are equal. An {\em $\alpha$-suborbit} is an orbit of
$G_{\alpha}$ on $\Omega$. If $G$ is transitive on $\Omega$, the {\em
subdegrees} of $G$ are the cardinalities of the $\alpha$-suborbits of
$G$ for some fixed $\alpha \in \Omega$.

There is a natural pairing between orbitals of $G$: if
$\Delta:=(\alpha, \beta)^G$ is an orbital, then its pair $\Delta^*$ is
the orbital $(\beta, \alpha)^G$. There is also a natural
correspondence between the orbital $\Delta$ and the $\alpha$-suborbit
$\Delta(\alpha)$, where $\Delta(\alpha):= \{\gamma \mid (\alpha,
\gamma) \in \Delta\}$. For every $\alpha$-suborbit $\Upsilon$ of $G$,
there is an orbital $\Delta$ such that $\Upsilon = \Delta(\alpha)$,
namely the orbital $\Delta:=(\alpha, \beta)^G$, where $\beta$ is some
vertex in $\Upsilon$. This correspondence is bijective. The pair of
the $\alpha$-suborbit $\Delta(\alpha)$ is the suborbit
$\Delta^*(\alpha)$; a suborbit or orbital is {\em self-paired} if it
is equal to its pair.\\

A {\it digraph} is a directed graph without multiple edges or loops;
it is a pair $(V\Gamma, A\Gamma)$, where $V\Gamma$ is the set of {\it
vertices} and $A\Gamma$ the set of {\it arcs} of $\Gamma$. The set
$A\Gamma$ consists of ordered pairs of distinct elements of $V\Gamma$.
Two vertices $\alpha, \beta \in V\Gamma$ are {\it adjacent} if either
$(\alpha, \beta)$ or $(\beta, \alpha)$ lies in $A\Gamma$. Throughout
this paper all paths will be undirected. The {\em distance} between
two connected vertices $\alpha$ and $\beta$ in $\Gamma$ is denoted by
$d_\Gamma(\alpha, \beta)$. A digraph is {\it locally-finite} if every
vertex is adjacent to at most finitely many other vertices.

If $G \leq \sym (\Omega)$ and $\Delta$ is an orbital of $G$, the
digraph $(\Omega, \Delta)$ is called an {\em orbital digraph of $G$}.
The connected components of an orbital digraph of $G$ are
$G$-invariant equivalence classes on $\Omega$; thus, if $G$ is
primitive, then every non-diagonal orbital digraph is connected.\\

Interest in the subdegree growth of infinite primitive permutation
groups possessing a finite suborbit whose pair is also finite stems
from the following observation.

\begin{thm} \label{thm:finite_suborbit_implies_all_suborbits_finite}
If $G$ is a group acting primitively on the infinite set $\Omega$,
possessing a finite suborbit whose pair is also finite, then every
suborbit of $G$ is finite and $\Omega$ is countable. \end{thm}

\begin{proof} Fix $\alpha \in \Omega$ and let $\beta^{G_\alpha}$ be a
finite $\alpha$-suborbit whose pair is finite. Let $\Gamma$ be the
orbital digraph $(\Omega, (\alpha, \beta)^G)$. This digraph is connected
since $G$ is primitive. Furthermore, the out-valency of $\alpha$ in
$\Gamma$ is equal to $|\beta^{G_\alpha}|$ and the in-valency of
$\alpha$ is equal to the size of the suborbit paired with
$\beta^{G_\alpha}$. Since both these sets are finite, the in-valency
and the out-valency of $\alpha$ are finite, so the valency of $\alpha$
is finite. Because $\Gamma$ is vertex-transitive, $\Gamma$ is locally
finite, so for all $r \geq 0$ the sphere $S_r(\alpha, \Gamma)$ is
finite, where $S_r(\alpha, \Gamma)$ denotes the set of vertices of
$\Gamma$ that lie at distance $r$ from $\alpha$.

Since $\Gamma$ is connected, $\Omega = \bigcup_{r=0}^\infty
S_r(\alpha, \Gamma)$. Thus $\Omega$ is a countable union of finite
sets, and is therefore countably infinite.

Since $G_\alpha$ fixes each sphere $S_r(\alpha, \Gamma)$ setwise,
every $\alpha$-suborbit of $G$ is finite. Because $G$ is primitive on
$\Omega$, it acts transitively on $\Omega$, so every suborbit of $G$
is finite.\end{proof}

Henceforth, any group whose suborbits are all finite will be described
as being {\it locally finite}. Thus a primitive group is locally
finite if and only if it has a locally finite orbital digraph; of
course, if it has a locally finite orbital digraph then all its orbital
digraphs will be locally finite.

Note this theorem is not true if we merely require a primitive group
$G$ to possess a finite suborbit. Indeed, in \cite{evans} Evans
constructs an infinite primitive group acting on an uncountably
infinite set $\Omega$, possessing a finite suborbit whose pair is
infinite.\\

The {\it complete digraph} on $t$ vertices, denoted by $K_t$, has $t$
vertices, between any two of which there is a pair of arcs, one in each direction.

The {\it connectivity} of a connected digraph $\Gamma$ is the
smallest possible size of a subset $W$ of $V\Gamma$ for which the
induced digraph $\Gamma \setminus W$ is disconnected. A {\it lobe} of
$\Gamma$ is a connected subgraph that is maximal subject to the
condition it has connectivity strictly greater than one. If $\Gamma$
has connectivity one, then the vertices $\alpha$ for which $\Gamma
\setminus \{\alpha\}$ is disconnected are called the {\it cut
vertices} of $\Gamma$.

A group acting on a digraph $\Gamma$ is said to be {\it
$k$-distance-transitive} if, for any four vertices $\alpha_1,
\alpha_2, \beta_1, \beta_2 \in V\Gamma$ with $d(\alpha_1, \alpha_2) =
d(\beta_1, \beta_2) \leq k$, there exists $g \in G$ such that
$\alpha_1^g = \beta_1$ and $\alpha_2^g = \beta_2$, and is {\em
distance-transitive} if it is $k$-distance-transitive for all positive
integers $k$. A digraph is {\it distance-transitive} if its automorphism
group acts upon it distance-transitively. The locally finite infinite
distance-transitive digraphs were classified by Macpherson
\cite{macpherson:dtg} and independently by Ivanov \cite{ivanov} thus.

\begin{thm} {\normalfont(\cite[Theorem 1.2]{macpherson:dtg})}
\label{thm:classification_of_dist_trans_graphs} An infinite locally finite
digraph $\Gamma$ is distance-transitive if and only if it is a regular
combinatorial tree, or $\Gamma$ has connectivity one, and for some integers
$m \geq 2$ and $t \geq 3$ the lobes of $\Gamma$ are isomorphic to the complete
digraph $K_t$, with each vertex in $\Gamma$ lying in $m$ lobes.\qed
\end{thm}

Adeleke and Neumann in \cite[Remark 29.8]{adeleke&neumann} observe
that if $G$ acts primitively on an infinite set $\Omega$ and has a
non-trivial self-paired finite suborbit of size $m$ then $\Omega$ is
countable, all the suborbits of $G$ are finite, $m \geq 2$, and if the
suborbits are arranged in a non-decreasing sequence $1 = m_0 \leq m_1
\leq \cdots$ then $m_r \leq (m-1)m_{r-1}$ for all sufficiently large
$r$. They then ask whether the subdegree growth rates of infinite
primitive groups that act distance-transitively on locally finite
distance-transitive digraphs are extremal and conjecture there might
exist a number $c$ which perhaps depends upon $G$, perhaps only on
$m$, such that $m_r \leq c(m-2)^{r-1}$.

In fact, this approach is naive, as it does not consider the existence
of a locally finite primitive group $G$, possessing at least two
distinct subdegrees $m$ and $m'$ such that infinitely many
$\alpha$-suborbits have cardinality $m$, and infinitely many have
cardinality $m'$. Such subdegrees will henceforth be said to occur
{\it infinitely often}. Indeed, if $G$ has at least two subdegrees,
each occurring infinitely often, of which $m$ is the smallest, then
under the above enumeration method $m_r = m$ for all sufficiently
large $r$. Any subdegree of $G$ that is strictly larger than $m$ would
therefore not be present in the subdegree sequence $(m_r)$.

In this paper we give examples of such groups, then define comprehensive
methods for enumerating the subdegrees of locally finite infinite primitive
groups and measuring their subdegree growth. Following this, we show that
certain rates of subdegree growth determine the structure of the group.

\section{Examples and constructions}

A {\em half-line} $L$ of a digraph $\Gamma$ is a one-way infinite
cycle-free path in $\Gamma$. The {\it ends} of $\Gamma$ are sets of half-lines, in which two
half-lines $L_1$ and $L_2$ lie in the same end if and only if there exist an infinite
number of pairwise-disjoint paths connecting a vertex in $L_1$ to a
vertex in $L_2$. In fact, this definition is an equivalence relation on the set of
half-lines of $\Gamma$, and the ends of $\Gamma$ are the equivalence classes of this  
relation.

It was noted in \cite{me:OrbGraphs} that if $G$ is an infinite primitive permutation group
possessing no infinite subdegree, then any two orbital digraphs of $G$ have the same ends.
Thus we define the {\it permutation-ends} of $G$ to be the set of ends of an orbital digraph
of $G$. It was also observed that all infinite locally finite primitive digraphs have one  
end,
or they have $2^{\aleph_0}$ ends; thus any infinite locally finite primitive permutation  
group must have precisely one permutation-end, or $2^{\aleph_0}$ permutation-ends.

In this section we give examples of locally finite infinite primitive groups with exactly  
one permutation-end, and examples with an
infinite number of permutation-ends. In both cases the groups in question posses at least two
distinct subdegrees, both of which occur infinitely often. 

It should be noted that the following examples may be used to construct several new examples  
of locally finite infinite primitive digraphs, which when taken to be undirected yield  
examples of locally finite infinite primitive graphs.

\subsection{Monster groups}
\label{subsection:monster}

In \cite{olshanski}, Ol'shanski\u\i \ shows for every prime $p >
10^{10}$, there is an infinite group in which all non-trivial proper
subgroups are of order $p$. We propose to call such groups
Tarski--Ol'shanski\u\i \ monster groups of order $p$. If $p$ is a
prime number greater than $10^{10}$ and $T_p$ is such a group, fix any
non-trivial proper subgroup $H \leq T_p$. Let $T_p$ act on the set of
right cosets $\Omega:=\{Hg \mid g \in T_p\}$ via
\[(Hg)^{g'} = Hgg'.\]
The kernel $K$ of this action is a normal subgroup of $T_p$, and is
therefore trivial. Indeed, suppose $K$ is a non-trivial, proper normal
subgroup of $T_p$. Then every non-trivial, proper subgroup of the
quotient group $T_p / K$ is of the form $H' / K$ for some proper
subgroup $H'$ of $T_p$ with $K < H'$. Now $T_p / K$ is infinite, and
every element has finite order, so $T_p / K$ contains a non-trivial,
proper subgroup $H' / K$, with $K < H' < T_p$. This is absurd,
however, since both $H'$ and $K$ must have order $p$. Hence $T_p$ is
simple, and its action on $\Omega$ is faithful.

The stabiliser of the coset $H.1 \in \Omega$ is the group $H \leq
T_p$. Since every subgroup of $T_p$ has order $p$, this group is a
maximal subgroup of $T_p$. Furthermore, every element of the finite
group $H$ has order $p$, so the orbits of $H$ acting on $\Omega
\setminus \{H\}$ must all have size $p$. Hence $T_p$ acts primitively
on $\Omega$, with all non-trivial suborbits finite of size $p$.

It was shown in \cite{me:OrbGraphs} that any infinite primitive
permutation group whose subdegrees are all finite that possesses an
orbital digraph with more than one end may be written as an amalgamated
free product. In particular, such groups must have at least one
element of infinite order. The group $T_p$ contains no elements of
infinite order, so every orbital digraph of this group must have precisely one
end. Thus, $T_p$ is an example of a locally finite primitive group
with one permutation-end, whose suborbits are all bounded above.

\subsection{Infinitely-ended constructions}
\label{section:constructions_with_infinitely_many_ends}

If $\Gamma$ is a digraph with connectivity one, with all lobes pairwise
isomorphic, and each vertex lies in precisely $m$ lobes, we shall denote
$\Gamma$ by $\Gamma(m, \Lambda)$, where $\Lambda$ is any lobe of $\Gamma$.

The bipartite digraph $T$ whose vertex set is the union of the set
$V_1$ of cut vertices of $\Gamma$ and the set $V_2$ of lobes of
$\Gamma$, in which $\alpha \in V_1$ is adjacent to $x \in V_2$ in $T$
if and only if $\alpha$ lies in $x$ in $\Gamma$, is in fact a tree, and
is known as the {\it block-cut-vertex tree} of $\Gamma$. Note that if
$\Gamma$ has connectivity one and block-cut-vertex tree $T$, then any
group $G$ acting on $\Gamma$ has a natural action on $T$.

Suppose we are given an integer $m \geq 2$, a locally finite primitive
arc-transitive digraph $\Lambda$ with connectivity at least two, and an
arc-transitive primitive non-regular group $H$ of automorphisms of
$\Lambda$. We will construct a primitive and arc-transitive group of
automorphisms $G$ of the digraph $\Gamma(m, \Lambda)$ such that the
subgroup of $\aut \Lambda$ induced by $G_{\{\Lambda\}}$ is equal to
$\overline{H}$, the closure of $H$ in $\aut \Lambda$ in the natural
complete topology on $\sym V\Gamma$. This construction will provide a
natural way of manufacturing large primitive groups from smaller ones.

A {\em relational structure} is a pair $(\Omega; \mathfrak{R})$, where
$\mathfrak{R}$ is a set of relations on the set $\Omega$. Suppose we
are given two relational structures $(\Omega_1; \mathfrak{R}_1)$ and
$(\Omega_2; \mathfrak{R}_2)$, and a bijective map $\phi:
\mathfrak{R}_1 \rightarrow \mathfrak{R}_2$. A map $\varphi : \Omega_1
\rightarrow \Omega_2$ is said to be an {\em isomorphism} between the
relational structures $(\Omega_1; \mathfrak{R}_1)$ and $(\Omega_2;
\mathfrak{R}_2)$ if and only if $\varphi(R) = \phi(R)$ for all
relations $R \in \mathfrak{R}_1$, where $\varphi(R)$ denotes the set
$\{(\varphi(\alpha_1), \ldots, \varphi(\alpha_n)) \mid (\alpha_1,
\ldots, \alpha_n) \in R\}$.

If $G$ is a permutation group acting on a non-empty set $\Omega$, and
$\Theta$ is an orbit of $G$ on $\Omega^n$, then it may be considered
to be a relation of arity $n$ on $\Omega$. Following
\cite{cameron:oligomorphic}, we define $\mathfrak{R}_G$ to be the set
of relations consisting of all orbits of $G$ on $\Omega^n$ where $n$
takes every value in the set of natural numbers. The {\em canonical
relational structure} associated with $G$ is the relational structure
$(\Omega; \mathfrak{R}_G)$.

\begin{thm} {\normalfont(\cite[Theorem 2.6]{cameron:oligomorphic})}
\label{thm:canonical_relational_structure} If $\Omega$ is a countably
infinite set then a group $G$ of permutations of $\Omega$ is closed if
and only if $G$ is the automorphism
group of the canonical relational structure of $G$. \qed \end{thm}

We begin by describing a relational structure that is based on the digraph
$\Gamma(m, \Lambda)$. Let $\Gamma$ denote this digraph, and, for any
digraph $\Gamma'$ with connectivity one, define $B(\Gamma')$ to be the
set of lobes of $\Gamma'$. Let $\mathfrak{R}$ be a set of relations on
the set $V\Lambda$ such that $A\Lambda \in \mathfrak{R}$. Since all
lobes of $\Gamma$ are isomorphic to $\Lambda$, we may, for each lobe
$\Delta$ of $\Gamma$, choose a digraph isomorphism
\[\varphi_\Delta : \Lambda \rightarrow \Delta.\]
For each relation $R \in \mathfrak{R}$ define $R_\Delta :=
\varphi_\Delta(R)$, where
\[\varphi_\Delta(R) := \{ ( \varphi_\Delta(\alpha_1), \ldots,
\varphi_\Delta(\alpha_n) ) \mid (\alpha_1, \ldots, \alpha_n) \in
R\}.\]
Now define $\mathfrak{R}_\Delta := \{R_\Delta \mid R \in
\mathfrak{R}\}$ and observe the relational structure $(V\Delta;
\mathfrak{R}_\Delta)$ is isomorphic to $(V\Lambda; \mathfrak{R})$.

For each vertex $\alpha \in V\Lambda$ the digraph $\Gamma \setminus
\{\alpha\}$ is not connected. One connected component of $\Gamma
\setminus \{\alpha\}$ contains all vertices in $\Lambda \setminus
\{\alpha\}$; the other connected components are disjoint from $\Lambda
\setminus \{\alpha\}$. Define $\Gamma_\alpha'$ to be the subgraph of
$\Gamma$ consisting of all connected components of $\Gamma \setminus
\{\alpha\}$ that are disjoint from $\Lambda \setminus \{\alpha\}$, and
let $\Gamma_\alpha$ be the subgraph of $\Gamma$ induced by
$V\Gamma_\alpha' \cup \{\alpha\}$. This digraph has connectivity one,
and for two vertices $\alpha, \beta \in V\Lambda$ the digraphs
$\Gamma_\alpha$ and $\Gamma_\beta$ are isomorphic.

Now fix $\alpha \in V\Lambda$ and let
\[\mathfrak{S}_\alpha := \bigcup_{\Delta \in B(\Gamma_\alpha)}
\mathfrak{R}_\Delta.\]
This is a set of relations on the set $V\Gamma_\alpha$, so we may
define $\Sigma_\alpha$ to be the relational structure
$(V\Gamma_\alpha; \mathfrak{S}_\alpha)$. For each $\gamma \in
V\Lambda$ choose a digraph isomorphism $\phi_\gamma : \Gamma_\alpha
\rightarrow \Gamma_\gamma$, with $\phi_\alpha$ the identity map. Now
define $\mathfrak{S}_\gamma:= \phi_\gamma(\mathfrak{S}_\alpha)$, where
$\phi_\gamma(\mathfrak{S}_\alpha) = \{\phi_\gamma(R) \mid R \in
\mathfrak{S}_\alpha\}$, and let $\Sigma_\gamma$ be the relational
structure $(V\Gamma_\gamma; \mathfrak{S}_\gamma)$. The mapping
$\phi_\gamma$ is thus an isomorphism between the relational structures
$\Sigma_\alpha$ and $\Sigma_\gamma$. This observation will underpin
the proof of our next lemma.

Now let
\[\mathfrak{S}:= \bigcup_{\gamma \in V\Lambda} \mathfrak{S}_\gamma,\]
and define $\Gamma(m, \Lambda, \mathfrak{R})$ to be the relational
structure $(V\Gamma; \mathfrak{S} \cup \mathfrak{R})$. This structure
has many important properties, the most useful being
\[\aut \Gamma(m, \Lambda, \mathfrak{R}) \leq \aut \Gamma(m, \Lambda).\]
For each lobe $\Delta$ of $\Gamma(m, \Lambda)$ the relational
structure $(V\Delta; \mathfrak{R}_\Delta)$ is isomorphic to
$(V\Lambda; \mathfrak{R})$. Henceforth, such relational structures
will be referred to as {\em relational lobes} of $\Gamma(m, \Lambda,
\mathfrak{R})$.

Recall that we have been given a primitive and non-regular group $H$
that acts arc-transitively on $\Lambda$. The canonical set of
relations $\mathfrak{R}_H$ consists of all orbits of $H$ on
$V\Lambda^n$, where $n$ takes every value in the set of natural
numbers. Let
\[G := \aut \Gamma(m, \Lambda, \mathfrak{R}_H).\]
We will show the subgroup of $\aut \Lambda$ induced by
$G_{\{\Lambda\}}$ is equal to $\overline{H}$, the closure of $H$ in
$\aut \Lambda$.

\begin{lemma} \label{lemma:infinitely_ended_construction} Let $T$ be
the block-cut-vertex tree of $\Gamma(m, \Lambda)$ and let $x$ be the
vertex of $T$ corresponding to $\Lambda$. If $\alpha \in V\Lambda$ and
$\beta \in V\Gamma$ lie in the same component of $T \setminus \{x\}$
then the subgroup of $\aut \Lambda$ induced by $G_{\alpha, \beta,
\{\Lambda\}}$ is isomorphic to $\overline{H}_\alpha$. \end{lemma}

\begin{proof} Let $G'$ be the subgroup of $\aut \Lambda$ induced by
$G_{\alpha, \beta, \{\Lambda\}}$. It is clear that $G' \leq \aut
(V\Lambda; \mathfrak{R}_H)$. By
Theorem~\ref{thm:canonical_relational_structure}, $\aut (V\Lambda;
\mathfrak{R}_H) = \overline{H}$. Thus,
\[G' \leq \overline{H}_\alpha.\]

Choose any automorphism $h \in \overline{H}_\alpha$. We will extend
this to an automorphism $\sigma$ of $\Gamma(m, \Lambda,
\mathfrak{R}_H)$ that fixes $\beta$.

For each pair of vertices $\gamma_1, \gamma_2 \in V\Lambda$ the
relational structures $\Sigma_{\gamma_1}$ and $\Sigma_{\gamma_2}$ are
isomorphic, so we may choose an isomorphism
\[\phi_{(\gamma_1, \gamma_2)} : \Sigma_{\gamma_1} \rightarrow
\Sigma_{\gamma_2},\]
which is equal to the identity map when $\gamma_1 = \gamma_2$. We now
construct a mapping $\sigma: \Gamma(m, \Lambda, \mathfrak{R}_H)
\rightarrow \Gamma(m, \Lambda, \mathfrak{R}_H)$ as follows. For each
$\delta \in V\Gamma$ there exists a unique vertex $\gamma \in
V\Lambda$ such that $\delta \in \Sigma_\gamma$, so set
\[\delta^\sigma := \delta^{\phi_{(\gamma, \gamma^h)}}.\]
It is simple to check this is a well-defined automorphism of
$\Gamma(m, \Lambda, \mathfrak{R}_H)$. Since $\alpha$ and $\beta$ lie
in the same component of $T \setminus \{x\}$, the vertex $\alpha$ must
be the unique element of $V\Lambda$ such that $\beta \in
\Sigma_\alpha$. Hence $\beta^\sigma = \beta$, so $G' =
\overline{H}_\alpha$.
\end{proof}

Before introducing our next result, a little notation is necessary.
Since there is a unique geodesic between any two vertices $\alpha$
and $\beta$ in a tree $T$, we shall denote the geodesic between
$\alpha$ and $\beta$ that includes both vertices by
$[\alpha, \beta]_T$; if we wish to exclude $\beta$ then we instead
write $[\alpha, \beta)_T$.

\begin{thm} \label{thm:infinitely_ended_construction} Suppose $H$ is a
primitive non-regular arc-transitive group of automorphisms of a
locally finite digraph $\Lambda$ with at least three vertices and
connectivity strictly greater than one. If $m \geq 2$ then the
automorphism group
$\aut \Gamma(m, \Lambda, \mathfrak{R}_H)$ acts primitively and
arc-transitively on the digraph $\Gamma(m, \Lambda)$. Furthermore, the
group induced by its action on $V\Lambda$ is $\overline{H}$, the
closure of $H$ in $\aut \Lambda$. \end{thm}

\begin{proof} We begin by showing the subgroup of $\aut \Lambda$
that is induced by
$\aut \Gamma(m, \Lambda, \mathfrak{R}_H)$ on $V\Lambda$ is
$\overline{H}$. Let $G = \aut \Gamma(m, \Lambda, \mathfrak{R}_H)$ and
let $G'$ be the subgroup of $\aut \Lambda$ induced by the action of
$G_{\{\Lambda\}}$ on $V\Lambda$. By construction, $G \leq \aut
\Gamma(m, \Lambda)$. Fix $\alpha \in V\Lambda$, and note that by
applying Lemma~\ref{lemma:infinitely_ended_construction} with $\beta =
\alpha$ we have $G'_\alpha = \overline{H}_\alpha$. Furthermore, $G'
\leq \aut(V\Lambda, \mathfrak{R}_H) = \overline{H}$, so
\[\overline{H}_\alpha \leq G' \leq \overline{H}.\]
Our choice of $\alpha$ was arbitrary, so we may choose a vertex
$\gamma \in V\Lambda$ that is distinct from $\alpha$ and note that
\[\overline{H}_\gamma \leq G' \leq \overline{H}.\]
Since $H$ acts primitively and non-regularly on $V\Lambda$, the group
$\overline{H}_\gamma$ does not fix $\alpha$; whence $G'$ does not fix
$\alpha$. Thus
\[\overline{H}_\alpha < G' \leq \overline{H}.\]
Since $H$ acts primitively on $V\Lambda$ the same must be true of
$\overline{H}$; therefore $\overline{H}_\alpha$ is a maximal subgroup
of $\overline{H}$. Consequently $G' = \overline{H}$.

It remains to show $G$ acts primitively and arc-transitively on the
digraph $\Gamma(m, \Lambda)$. By construction $G$ acts transitively on
the relational lobes of $\Gamma(m, \Lambda, \mathfrak{R}_H)$, and
therefore acts transitively on the set of lobes of the digraph
$\Gamma(m, \Lambda)$.

The action of $G$ on the vertices of $\Gamma$ is transitive. Indeed,
suppose $\delta$ is any vertex in $V\Gamma$. Then there exists a lobe
$\Delta$ of $\Gamma$ containing $\delta$. Choose $g \in G$ such that
$\Delta^g = \Lambda$. Since $G' = \overline{H}$, and $H$ acts
transitively on the vertices of $\Lambda$, there exists an
automorphism $h \in G_{\{\Lambda\}}$ such that $\delta^{gh} = \alpha$.

Furthermore, the stabiliser $G_\alpha$ transitively permutes the lobes
of $\Gamma$ that contain $\alpha$. For, suppose $\Lambda'$ is a lobe
of $\Gamma$ containing $\alpha$. Then there is an automorphism $g \in
G$ such that $\Lambda^g = \Lambda'$. We again observe that since $G' =
\overline{H}$, and $H$ acts transitively on the vertices of $\Lambda$,
there exists an automorphism $h \in G_{\{\Lambda\}}$ such that
$\alpha^h = \alpha^{g^{-1}}$. Thus $\Lambda^{h g} = \Lambda'$ and
$\alpha^{h g} = \alpha$.

Since $H$ acts vertex- and arc-transitively on $\Lambda$, and
$G_\alpha$ transitively permutes the lobes of $\Gamma$ that contain
$\alpha$, we must have $G$ is arc-transitive on $\Gamma$.

Finally, let $\rho$ be a non-trivial $G$-congruence on $V\Gamma$. Let
$T$ be the block-cut-vertex tree of $\Gamma$. Choose $\beta \in
\rho(\alpha) \setminus \{\alpha\}$ of minimal distance in $T$ from
$\alpha$. We claim $\alpha$ and $\beta$ lie in a common lobe of
$\Gamma$. Indeed, suppose this is not the case. Let $x$ be the vertex
adjacent to $\beta$ in the geodesic $[\alpha, \beta]_T$ between $\alpha$
and $\beta$, and let
$\gamma$ be the vertex adjacent to $x$ in $[\alpha, \beta)_T$. The
vertex $x$ corresponds to a lobe $\Lambda'$ of $\Gamma$, which, since
$\alpha$ and $\beta$ lie in distinct lobes, is not equal to $\Lambda$.
Because $G$ acts transitively on the lobes of $\Gamma$, the group
induced on $V\Lambda'$ by $G_{\alpha, \gamma, \{\Lambda'\}}$ is
isomorphic to $\overline{H}_\alpha$ acting on $V\Lambda$. Since
$\overline{H}$ is primitive and non-regular on $V\Lambda$, there
exists $g \in G_{\alpha, \gamma, \{\Lambda'\}}$ such that $\beta^g
\not = \beta$. Thus $\beta, \beta^g \in \rho(\alpha)$ share a common
lobe, and our claim is established. So, without loss of generality, we
may assume $\alpha$ and $\beta$ lie in $\Lambda$. Now $\rho$ induces a
$\overline{H}$-congruence on $V\Lambda$, so we must have $V\Lambda
\subseteq \rho(\alpha)$. Furthermore, $G_\alpha$ acts transitively on
the lobes of $\Gamma$ that contain $\alpha$, so $\rho$ is the
universal relation. Hence $G$ acts primitively on $\Gamma$.
\end{proof}

If $H$ is a primitive and non-regular group of permutations of a set
$\Omega$, with $|\Omega| \geq 3$, possessing an orbital digraph
$\Lambda$ with connectivity strictly greater than one, then the
primitive group $G$ constructed above will be called the {\em $m$-fold
graph product of $H$}, and will be denoted by $G(m, H)$. By the above
theorem, this group acts primitively on the vertex set of $\Gamma(m,
\Lambda)$. Since this digraph has infinitely many ends, and is an
orbital digraph of $G(m, \Lambda)$, all orbital digraphs of $G$ will have
infinitely many ends.

\begin{example} \label{example:infintely_ended} Let $p$ be a prime
number with $p > 10^{10}$ and let $T_p$ be a Tarski--Ol'shanski\u\i \
monster group of order $p$. Recall this group acts primitively on its
coset space $\Omega$. Let $\Lambda$ be an orbital digraph of $T_p$. Then
$\Lambda$ is a one-ended primitive digraph, and therefore has
connectivity strictly greater than one.

By Theorem~\ref{thm:infinitely_ended_construction}, the group $G(m,
T_p)$ acts primitively on $\Gamma(m, \Lambda)$, and has infinitely
many distinct subdegrees, each occurring infinitely often. In fact,
for all $\alpha \in V\Gamma$, the $\alpha$-subdegrees of $G(m, T_p)$
are $m (m-1)^{r-1} p^r$ for $r \geq 1$. For each positive integer $r$
there are infinitely many $\alpha$-suborbits of $G(m, T_p)$ with
cardinality $m (m-1)^{r-1} p^r$. \end{example}

\subsection{One-ended constructions}

The following result allows one to construct infinitely many examples
of locally finite primitive groups with one permutation-end.

\begin{thm} \label{thm:wreath_product_construction} If $G$ is a
locally finite primitive group of permutations of an infinite set
$\Omega$ and $m \geq 2$, then the wreath product $G \Wr \sym(m)$ under
the product action is a locally finite primitive group of permutations
of $\Omega^m$ with one permutation-end. \end{thm}

This result may be deduced from the following lemmas. Fix a group $G$
acting primitively on an infinite set $\Omega$, possessing a finite
suborbit whose pair is also finite. Let $\Gamma$ be an orbital digraph
of $G$ on $\Omega$, and define $H:=G \Wr \sym(m)$. It is well known
that $H$ acts primitively on the set $\Omega^m$ under the product
action; see, for example \cite[Lemma 2.7A]{dixon&mortimer}. Fix
$\alpha, \beta \in \Omega$ such that $\alpha$ and $\beta$ are adjacent
in $\Gamma$, with $(\alpha, \beta) \in A\Gamma$. Let
$\underline{\alpha}:= (\alpha, \ldots, \alpha) \in \Omega^m$ and
$\underline{\beta}:=(\beta, \alpha, \ldots, \alpha) \in \Omega^m$, and
define
\[\Lambda := (\Omega^m, (\underline{\alpha}, \underline{\beta})^H),\]
where $H$ acts on $\Omega^m$ via the product action.

\begin{lemma} \label{lemma:S_1_of_Gamma_Wr_Sm} The vertex $(\gamma_1,
\ldots, \gamma_m)$ is adjacent to $\underline{\alpha}$ in $\Lambda$ if
and only if \[\sum_{i=1}^m d_{\Gamma}(\alpha, \gamma_i) = 1.\]
\end{lemma}

\begin{proof} Let $\underline{\gamma} := (\gamma_1, \ldots, \gamma_m)
\in V\Lambda$. Since $\Lambda$ is an orbital digraph of $H$ on
$\Omega^m$, it is arc transitive. Therefore $\underline{\gamma}$ lies
in the sphere $S_1(\underline{\alpha}, \Lambda)$ if and only if there
exists $g = (g_1, \ldots, g_m)\sigma \in H$ such that
$(\underline{\alpha}, \underline{\beta})^g$ is equal to
$(\underline{\alpha}, \underline{\gamma})$ or $(\underline{\gamma},
\underline{\alpha})$. If $(\underline{\alpha}, \underline{\beta})^g =
(\underline{\alpha}, \underline{\gamma})$ then $g \in
H_{\underline{\alpha}}$ and so $g_1, \ldots, g_m \in G_{\alpha}$. Thus
$\underline{\gamma} = \underline{\beta}^g = (\beta^{g_1}, \alpha,
\ldots, \alpha)^\sigma$, and therefore $\sum_{i=1}^m
d_{\Gamma}(\alpha, \gamma_i) = d_{\Gamma}(\alpha, \beta^{g_1}) = 1$.
Otherwise, if $(\underline{\alpha}, \underline{\beta})^g =
(\underline{\gamma}, \underline{\alpha})$ then $\underline{\beta}^g =
\underline{\alpha}$, so $g_2, \ldots, g_m \in G_{\alpha}$ but
$\beta^{g_1} = \alpha$. Thus $\underline{\gamma} =
\underline{\alpha}^g = (\alpha^{g_1}, \alpha, \ldots, \alpha)^\sigma$
and $\sum_{i=1}^m d_{\Gamma}(\alpha, \gamma_i) = d_{\Gamma}(\alpha,
\alpha^{g_1}) = 1$.

Conversely, suppose $\underline{\gamma} = (\gamma_1, \ldots, \gamma_m)
\in V\Lambda$ and $\sum_{i=1}^m d_{\Gamma}(\alpha, \gamma_i) = 1$.
Then there exists a unique value of $i$ such that $d_{\Gamma}(\alpha,
\gamma_i) = 1$ and $\gamma_j = \alpha$ for all $j \not = i$. Thus,
there exists $g_1 \in G$ such that $(\alpha, \beta)^{g_1}$ is equal to
$(\alpha, \gamma_i)$ or $(\gamma_i, \alpha)$. Let $\sigma \in \sym(m)$
be the $2$-cycle $(1 i)$ and set $g := (g_1, 1, \ldots, 1)\sigma \in
H$. If $(\alpha, \beta)^{g_1} = (\alpha, \gamma_i)$ then $g \in
H_{\underline{\alpha}}$ and $\underline{\beta}^g = (\beta^{g_i},
\alpha, \ldots, \alpha)^\sigma = (\gamma_i, \alpha, \ldots,
\alpha)^\sigma = \underline{\gamma}$, and therefore
$d_{\Lambda}(\underline{\alpha}, \underline{\gamma}) = 1$. On the
other hand, if $(\alpha, \beta)^{g_1} = (\gamma_i, \alpha)$ then
$\underline{\beta}^g = \underline{\alpha}$ and $\underline{\alpha}^g =
(\alpha^{g_1}, \alpha, \ldots, \alpha)^\sigma = (\gamma_i, \alpha,
\ldots, \alpha)^\sigma = \underline{\gamma}$, so again
$d_{\Lambda}(\underline{\alpha}, \underline{\gamma}) = 1$. We have
thus shown $(\gamma_1, \ldots, \gamma_m) \in S_1((\alpha, \ldots,
\alpha), \Lambda)$ if and only if $\sum_{i=1}^m d_{\Gamma}(\alpha,
\gamma_i) = 1$. Since $\Lambda$ is vertex-transitive the hypothesis
holds for any vertex $(\delta_1, \ldots, \delta_m) \in V\Lambda$.
\end{proof}

\begin{lemma} \label{lemma:wreath_prod_graph_sphere} The vertex
$(\gamma_1, \ldots, \gamma_m)$ lies in $S_r((\delta_1,
\ldots, \delta_m), \Lambda)$ if and only if
\[\sum_{i=1}^m d_{\Gamma}(\delta_i, \gamma_i) = r.\] \end{lemma}

\begin{proof} We proceed by induction. Since $H$ is vertex-transitive,
the hypothesis holds when $r=1$ by
Lemma~\ref{lemma:S_1_of_Gamma_Wr_Sm}. Fix $\underline{\delta} =
(\delta_1, \ldots, \delta_m) \in V\Lambda$ and $k > 1$ and suppose the
hypothesis is true for all $r \leq k$. Choose $\underline{\gamma} =
(\gamma_1, \ldots, \gamma_m) \in V\Lambda$ with $\sum_{i=1}^m
d_{\Gamma}(\delta_i, \gamma_i) = k+1$. Since $\sum_{i=1}^m
d_{\Gamma}(\delta_i, \gamma_i) > 2$ there exists $j$ for which one may
choose $\gamma_j^\prime \in S_1(\gamma_j, \Gamma)$ such that
$d_{\Gamma}(\delta_j, \gamma_j^\prime) = d_{\Gamma}(\delta_j,
\gamma_j) - 1$. Let $\gamma_i^\prime:=\gamma_i$ for all $i \not = j$
and put $\underline{\gamma^\prime}:=(\gamma_1^\prime, \ldots,
\gamma_m^\prime)$. Now $\sum_{i=1}^m d_{\Gamma}(\delta_i,
\gamma_i^\prime) = k$, so by assumption, $\underline{\gamma^\prime}
\in S_k(\underline{\delta}, \Lambda)$; furthermore,
$d_{\Lambda}(\underline{\gamma}, \underline{\gamma^\prime}) = 1$, so
$\underline{\gamma} \in S_{k-1}(\underline{\delta}, \Lambda) \cup
S_k(\underline{\delta}, \Lambda) \cup S_{k+1}(\underline{\delta},
\Lambda)$. If $\underline{\gamma}$ lies in
$S_{k-1}(\underline{\delta}, \Lambda)$ or $S_k(\underline{\delta},
\Lambda)$ then by assumption $\sum_{i=1}^m d_{\Gamma}(\delta_i,
\gamma_i)$ is equal to $k-1$ or $k$ respectively; since this is not
the case, we must have $\underline{\gamma} \in
S_{k+1}(\underline{\delta}, \Lambda)$.

Conversely, suppose $\underline{\gamma} = (\gamma_1, \ldots, \gamma_m)
\in S_{k+1}(\underline{\delta}, \Lambda)$. Since
$d_{\Lambda}(\underline{\delta}, \underline{\gamma}) > k$ we have
$\sum_{i=1}^m d_{\Gamma}(\delta_i, \gamma_i) \geq k+1$ by the
induction hypothesis. Since $\Lambda$ is connected, there exists
$\underline{\gamma^\prime} = (\gamma_1^\prime, \ldots,
\gamma_m^\prime) \in S_k(\underline{\delta}, \Lambda) \cap
S_1(\underline{\gamma}, \Lambda)$. Now $\sum_{i=1}^m
d_{\Gamma}(\delta_i, \gamma_i^\prime) = k$ and $\sum_{i=1}^m
d_{\Gamma}(\gamma_i, \gamma_i^\prime) = 1$. For each $i$ we have
$d_{\Gamma}(\delta_i, \gamma_i) \leq d_{\Gamma}(\delta_i,
\gamma_i^\prime) + d_{\Gamma}(\gamma_i^\prime, \gamma_i)$. Whence,
$\sum_{i=1}^m d_{\Gamma}(\delta_i, \gamma_i) \leq \sum_{i=1}^m
(d_{\Gamma}(\delta_i, \gamma_i^\prime) + d_{\Gamma}(\gamma_i^\prime,
\gamma_i)) = k + 1$. Hence $\sum_{i=1}^m d_{\Gamma}(\delta_i,
\gamma_i) = k+1$. \end{proof}

\begin{lemma} The digraph $\Lambda$ is infinite, primitive, locally
finite and arc-transitive, with one end.\end{lemma}

\begin{proof} The group $H$ acts primitively on $V\Lambda$, so the
digraph $\Lambda$ is primitive. It is infinite because $V\Lambda =
\Omega^m$ is infinite, and is locally finite by
Lemma~\ref{lemma:S_1_of_Gamma_Wr_Sm}; furthermore, since $\Lambda$ is
an orbital digraph of $H$ on $\Omega^m$, it is arc-transitive.

It remains to prove that $\Lambda$ has one end. We will show, for all
$r \geq 1$, given any pair of vertices $\underline{\gamma},
\underline{\delta} \in S_{r+1}(\underline{\alpha}, \Lambda)$, there is
a path connecting $\underline{\gamma}$ to $\underline{\delta}$ that is
not contained in the ball $B_r(\underline{\alpha}, \Lambda)$, where
$B_r(\underline{\alpha}, \Lambda)$ is the set of all vertices in
$\Lambda$ whose distance from $\alpha$ is at most $r$. From
this, we may deduce there is no finite subgraph of $\Lambda$ that one
may remove to leave at least two disjoint infinite connected
components. Whence, $\Lambda$ has precisely one end.

We begin by observing that, given any two vertices $(\gamma_1, \ldots,
\gamma_m)$ and $(\delta_1, \ldots, \delta_m)$ in $V\Lambda$, for any
path in $\Gamma$ between $\gamma_1$ and $\delta_1$ there exists a
corresponding path in $\Lambda$ between $(\gamma_1, \gamma_2, \ldots,
\gamma_m)$ and $(\delta_1, \gamma_2, \ldots, \gamma_m)$. Indeed, by
Lemma~\ref{lemma:S_1_of_Gamma_Wr_Sm}, if $\xi$ lies on the path in
$\Gamma$ between $\gamma_1$ and $\delta_1$, then $(\xi, \gamma_2,
\ldots, \gamma_m)$ lies on the corresponding path in $\Lambda$ between
$(\gamma_1, \gamma_2, \ldots, \gamma_m)$ and $(\delta_1, \gamma_2,
\ldots, \gamma_m)$. This observation can also be made for paths
between $\gamma_i$ and $\delta_i$ for all $i$ satisfying $1 \leq i
\leq m$.

Fix $r \geq 1$ and two distinct vertices $\underline{\gamma}  =
(\gamma_1, \ldots, \gamma_m)$ and $\underline{\delta} = (\delta_1,
\ldots, \delta_m)$ in $S_{r+1}(\underline{\alpha}, \Lambda)$. Let
$B_r:= B_r(\underline{\alpha}, \Lambda)$. We will describe four
vertices $\underline{\xi_1}, \ldots, \underline{\xi_4} \in V\Lambda
\setminus B_r$ such that there exist paths in $\Lambda$ between
$\underline{\gamma}$ and $\underline{\xi_1}$; between
$\underline{\delta}$ and $\underline{\xi_4}$; and between
$\underline{\xi_i}$ and $\underline{\xi_{i+1}}$ for $1 \leq i < 4$
that are all disjoint from $B_r$, thus showing there exists a path in
$\Lambda \setminus B_r$ from $\underline{\gamma}$ to
$\underline{\delta}$.

Let $d_1:= d_\Gamma(\alpha_1, \delta_1) - 1$ and $d_m:=
d_\Gamma(\alpha_m, \gamma_m) - 1$. Choose a vertex $\gamma_m' \in
V\Gamma \setminus B_r(\alpha_m, \Gamma)$ such that there is a path in
$\Gamma$ between $\gamma_m$ and $\gamma_m'$ that is disjoint from
$B_{d_m}(\alpha_m, \Gamma)$. Similarly, choose a vertex $\delta_1' \in
V\Gamma \setminus B_r(\alpha_1, \Gamma)$ such that there is a path in
$\Gamma$ between $\delta_1$ and $\delta_1'$ that is disjoint from
$B_{d_1}(\alpha_1, \Gamma)$. Now define
\begin{align*} \underline{\xi_1} &:= (\gamma_1, \gamma_2, \ldots,
\gamma_{m-1}, \gamma_m'); \\
\underline{\xi_2} &:= (\delta_1, \delta_2, \ldots, \delta_{m-1}, \gamma_m'); \\
\underline{\xi_3} &:= (\delta_1', \delta_2, \ldots, \delta_{m-1}, \gamma_m'); \\
\underline{\xi_4} &:= (\delta_1', \delta_2, \ldots, \delta_{m-1},
\delta_m). \end{align*}

The path in $\Gamma$ between $\gamma_m$ and $\gamma_m'$ that is
disjoint from $B_{d_m}(\alpha_m, \Gamma)$ corresponds to a path in
$\Lambda$ between the vertices $(\gamma_1, \ldots, \gamma_{m-1}, \gamma_m)$ and
$(\gamma_1, \ldots, \gamma_{m-1}, \gamma_m')$ which, by
Lemma~\ref{lemma:wreath_prod_graph_sphere}, is disjoint from $B_r$.
Hence, there exists a path in $\Lambda \setminus B_r$ between
$\underline{\gamma}$ and $\underline{\xi_1}$. A similar argument shows
 there exists a path in $\Lambda \setminus B_r$ between
$\underline{\delta}$ and $\underline{\xi_4}$.

Finally, observe that any vertex $\underline{\chi} = (\chi_1, \ldots,
\chi_m) \in V\Lambda$ satisfying $\chi_i \in V\Gamma \setminus
B_r(\alpha_i, \Gamma)$ for some $i$ does not lie in $B_r$ by
Lemma~\ref{lemma:wreath_prod_graph_sphere}. Therefore, there exist
paths in $\Lambda \setminus B_r$ between $\underline{\xi_1}$ and
$\underline{\xi_2}$; between $\underline{\xi_2}$ and
$\underline{\xi_3}$; and between $\underline{\xi_3}$ and
$\underline{\xi_4}$. Hence, $\Lambda$ has precisely one end.
\end{proof}

Thus the proof of Theorem~\ref{thm:wreath_product_construction} is complete.

\begin{example} \label{example:height_2} Let $p$ be a prime number
greater than $ 10^{10}$ and let $T_p$ be a Tarski--Ol'shanski\u\i \
monster group of order $p$. Recall this group acts primitively on its
coset space $\Omega$. Let $G:= T_p \Wr \sym (2)$. This group acts
primitively on the set $\Omega^2$ by
Theorem~\ref{thm:wreath_product_construction}, and has precisely three
non-trivial subdegrees, $2p, p^2$ and $2 p^2$, each occurring
infinitely often. \end{example}

\begin{example} The orbital digraphs of the group $G(m, T_p) \Wr
\sym(2)$ have precisely one end. Furthermore, this group has
infinitely many distinct subdegrees, each occurring infinitely often.
\end{example}

\section{Enumeration and growth}

\subsection{The upper and lower subdegree sequences}
\label{section:upper_and_lower_growth}

As the above examples show, the concept of subdegree growth requires
re-examining. We begin by clarifying some terminology. Let $G$ act
primitively on an infinite set $\Omega$, possessing a finite suborbit
whose pair is also finite. Then every suborbit of $G$ is finite. Fix
$\alpha \in \Omega$.

The {\em set of subdegrees} of $G$ is defined to be the set whose
elements are the cardinalities of its $\alpha$-suborbits. The {\em
multiset of subdegrees} of $G$ is a function $\mu$ from the set of
subdegrees to the extended non-negative integers defined as follows.
If $m$ is any subdegree of $G$ then $\mu(m)$ is the number of
$\alpha$-suborbits with cardinality $m$. Having defined this multiset
formally, we shall speak of it informally as a set of elements of the
set of subdegrees, in which some elements occur with multiplicity
greater than one.

For our given group $G$, there exists a minimal ordinal number
$\kappa$ such that one may enumerate all elements of the multiset of
subdegrees as a monotonic increasing sequence $(m_\gamma)$ for $\gamma
< \kappa$. This sequence $(m_\gamma)$ is called the {\em subdegree
sequence of $G$}, and the ordinal number $\kappa$ is called the {\em
height} of $G$. The precise definition of the subdegree sequence of
$G$ is now hopefully clear: $(m_\gamma)$ is a non-decreasing sequence
of elements of the set of subdegrees, in which each subdegree $m$
appears in the sequence with multiplicity $\mu(m)$.

By Theorem~\ref{thm:finite_suborbit_implies_all_suborbits_finite}, the
set of subdegrees of $G$ is a finite or countably infinite set. Since
each entry occurs with multiplicity one, we may enumerate all its
elements $M_1 < M_2 < \cdots$. The sequence $(M_r)$ is called the {\em
upper subdegree sequence} of $G$. The {\em lower subdegree sequence}
of $G$ is the sequence $(m_r)_{r < \omega}$. Note that both the upper
and lower subdegree sequences are subsequences of the subdegree
sequence of $G$, and are indexed by the natural numbers.

\begin{lemma} If $G$ is a locally finite primitive group with height $h$ then
\[\omega \leq h \leq \omega^2.\] \end{lemma}

\begin{proof} Let $X$ and $Y$ be the set and multiset of subdegrees of
$G$ respectively. Since $G$ acts on an infinite set, and all
subdegrees of $G$ are finite, $Y$ is infinite, so the height of $G$ is
at least $\omega$.

Each subdegree in $X$ occurs with multiplicity at most $\aleph_0$ in
$Y$, so the subdegree sequence can be enumerated in a non-decreasing
sequence of length at most $\omega^2$. \end{proof}

Subdegree growth is similar in many ways to the growth of connected
locally finite infinite digraphs, an area of research that is already
well-establish. Following \cite{watkins:survey}, if $t_1 \leq t_2 \leq
t_3 \leq \cdots$ is a sequence of positive real numbers, we define the
concept of {\em growth} as follows. If there exist positive real
numbers $c_1$, $c_2$ and $d \geq 1$ such that
\[c_1 r^d \leq t_r \leq c_2 r^d,\]
we say that the sequence has {\it polynomial growth of degree $d$}. A
sequence exhibiting polynomial growth of degree $0$ is often called
{\em bounded}. The growth is {\em subexponential} if, for all $a > 1$,
\[\liminf_{r \rightarrow \infty} \frac{t_r}{a^r}=0;\]
it is {\em exponential} if there exists a constant $a>1$ such that
$\displaystyle{\liminf_{r \rightarrow \infty} t_r/a^r}$ and
$\displaystyle{\limsup_{r \rightarrow \infty} t_r/a^r}$ are non-zero
and finite. The growth is said to be {\em super exponential} if
$\displaystyle{\liminf_{r \rightarrow \infty} t_r/a^r}$ is infinite
for all $a > 1$.

The {\em growth} of an infinite locally finite vertex-transitive
connected digraph $\Gamma$ is the growth of the sequence $(|B_r(\alpha,
\Gamma)|)$. A great deal of work has been done in this area. It was
shown in \cite{imrich&seifter:linear_growth} that it is precisely the
two-ended digraphs that have growth of degree one. Furthermore, Trofimov
has shown in \cite{trofimov:polynomial_growth} that if a locally
finite digraph has polynomial growth then it is not primitive. All
locally finite infinitely-ended vertex-transitive digraphs have
exponential sphere growth, while it is noted in \cite{watkins:survey}
that one-ended digraphs exhibit all possible rates of growth of degree
at least $2$, except super exponential.

By examining the growth of both the lower and upper subdegree
sequences, one may obtain similar results pertaining to subdegree
growth rates. It is natural to consider the subdegree growth of
primitive groups in this way. Indeed, the intuitive but flawed
approach taken by Adeleke and Neumann in \cite[Remark
29.8]{adeleke&neumann} only considered the existence of primitive
groups with height $\omega$. For such groups, the lower subdegree
sequence is equal to the subdegree sequence (although it is not
necessarily equal to the upper subdegree sequence), so the two
approaches are equivalent in this case.

For groups with height $\omega$, the growth of the upper subdegree
sequence is of secondary importance, as all relevant subdegree
information can be found in the lower subdegree sequence. However,
this is not the case for groups with height strictly greater than
$\omega$.

\begin{lemma} If an infinite primitive group $G$ whose subdegrees are
all finite does not have height $\omega$, then the lower subdegree
sequence is always bounded. \end{lemma}

\begin{proof} Suppose $G$ has height $h > \omega$, and let
$(m_\gamma)$ be the subdegree sequence of $G$. Then there are finite
constants $c_1 := 1$ and $c_2 := m_{\omega}$ such that $c_1 \leq m_r
\leq c_2$ for all $r < \omega$. \end{proof}

The growth of the upper subdegree sequence will be used only to
differentiate between the subdegree growth of groups exhibiting a
bounded lower subdegree sequence; thus, when referring to the
subdegree growth of a group, unless otherwise stated, it is the lower
subdegree growth to which we are referring.

The constructions detailed in this chapter may be used to create
myriad examples of groups with exponential, subexponential and
polynomial growth, several of which are given below. The final example
illustrates how they may also be used to show that the above list of
possible rates of growth is not exhaustive.

\begin{example} Fix $m \geq 2$ and $t \geq 3$, and let $G:=\aut
\Gamma(m, K_{t+1})$. The group $G$ acts distance-transitively on
$V\Gamma$, and has height $\omega$. The upper and lower subdegree
sequences of $G$ are equal, with $M_r = m_r = m (m-1)^{r-1}t^r$, for
all $r \in \mathbb{N}$.\end{example}

\begin{example} If $p$ is prime, and $p > 10^{10}$, then the group
$G(m, T_p)$ has a bounded lower subdegree sequence, with $m_r = m p$
for all integers $r \geq 1$. However, $M_r = m (m-1)^{r-1}p^r$ for all
$r \geq 1$. Intuitively, while this group exhibits slow subdegree
growth when compared to groups with non-bounded lower subdegree
sequences, when instead it is compared with other groups possessing
bounded lower subdegree sequences its subdegree growth is extremely
fast. \end{example}

\begin{example} Given a prime integer $p$ with $p > 10^{10}$, both
$T_p$ and the wreath product $T_p \Wr \sym (2)$ have bounded upper and
lower subdegree sequences. The subdegree growth of these groups is
therefore, intuitively, very slow. Indeed, $T_p$ provides a lower
bound on both upper and lower subdegree growth rates. \end{example}

\begin{example} \label{ex:nonPolyNonExp} Fix integers $m, t \geq 2$
and define $\Gamma := \Gamma(m, K_{t+1})$. Let $G:= \aut \Gamma$ and
let $H:= G \Wr \sym(2)$. This group acts primitively on the set
$\Omega := V\Gamma \times V\Gamma$, with all subdegrees finite. Choose
vertices $\alpha$ and $\beta$ adjacent in $\Gamma$ and let $\Lambda$
be the orbital digraph $(\Omega, ((\alpha, \alpha),(\alpha, \beta))^H)$.
We shall denote the number of $H_{(\alpha, \alpha)}$-orbits in the
sphere $S_r((\alpha, \alpha), \Lambda)$ by $n_r$, and the number of
$H_{(\alpha, \alpha)}$-orbits in the ball $B_r((\alpha, \alpha),
\Lambda)$ by $N_r$.

Since $G$ acts distance-transitively on $\Gamma$, $n_r = \lceil r/2 \rceil$,
where $\lceil r/2 \rceil$ denotes the smallest integer greater than or
equal to $r/2$. Furthermore, the subdegrees of $H$ in the sphere
$S_r((\alpha, \alpha), \Lambda)$ are $2m(m-1)^{r-1}t^r$ and $2 m^2
(m-1)^{r-2}t^r$ for $r \geq 2$, and, if $r$ is even, $m^2(m-1)^{r -
2}t^{r}$.

The largest suborbit in the sphere $S_r((\alpha, \alpha), \Lambda)$
has size $2m^2(m-1)^{r-2}t^r$, and for all sufficiently large integers
$r$,
\[2m^2(m-1)^{r-2}t^r \leq (m-1)^{2r-2}t^{2r}.\]
Thus, there exists an integer $R$ such that, for all $r > R$, the
largest suborbit in $S_r((\alpha, \alpha), \Lambda)$ has cardinality
strictly less than the cardinality of every suborbit in $\Lambda
\setminus B_{2r}((\alpha, \alpha), \Lambda)$.

Now consider the lower subdegree sequence $(m_r)$ of H which, in this
case, is equal to the subdegree sequence of $H$. Choose $r > R$ and
find the largest integer $s \geq 2$ such that
\[m_s = 2m^2(m-1)^{r-2}t^r.\]
Our aim is to find the number of subdegrees that are less than or
equal to $m_s$, and from this determine $s$.

Since $m_s$ is the largest suborbit in $B_r((\alpha, \alpha),
\Lambda)$ we have $s \geq N_r$. Furthermore, $m_s$ is strictly less
than the cardinality of every $(\alpha, \alpha)$-suborbit in $\Lambda
\setminus B_{2r}((\alpha, \alpha), \Lambda)$, so there can be at most
$N_{2r}$ suborbits with cardinality less than or equal to $m_s$. Hence
$s \leq N_{2r}$, and therefore $N_r \leq s \leq N_{2r}$. Since $N_r =
\sum_{i=1}^r n_i$, we have $N_r \geq r(r+1)/4$ and $N_r \leq
r(r+1)/2$. Thus
\[r(r+1)/4 \leq s \leq r(2r+1).\]

If the lower subdegree sequence of $H$ exhibits polynomial growth of
degree $d \geq 1$ then there exist positive real numbers $c_1, c_2$
such that for all $k \geq 1$,
\[c_1 k^d \leq m_k \leq c_2 k^d.\]
However, for all $r>R$ we may choose a maximal integer $s$ with $m_s =
2m^2(m-1)^{r-2}t^r$. For all sufficiently large $r$ we therefore have
$s^d \leq r^d(2r+1)^d \leq 2m^{2d}(m-1)^{d(r-2)}t^{dr}$, so the lower
subdegree growth rate of $H$ is not polynomial.

Although the growth of the lower subdegree sequence of $H$ is faster
than polynomial growth, it is not exponential. Indeed, for all $r > R$
there exists an integer $s_r$ with $m_{s_r} = 2 m^2 (m-1)^{r-2}t^r$
and $s_r \geq r(r+1)/4$. Hence, for all $a > 1$,
\[ \lim_{r \rightarrow \infty} \frac{m_{s_r}}{a^{s_r}} \leq
\frac{2m^2(m-1)^{r-2}t^r}{a^{r(r+1)/4}} = 0.\]

The group $H$ is thus an example of a group exhibiting subexponential,
non-polynomial growth. The existence of such a group demonstrates that
the list of possible growth rates given previously is not exhaustive.
\end{example}

\subsection{The average subdegree sequence}

At this point it is perhaps relevant to draw the reader's attention to
another seemingly natural method for enumerating subdegrees, that
neatly avoids the problems caused by subdegrees which occur infinitely
often. Let $G$ act primitively on an infinite set $\Omega$, and
suppose every subdegree of $G$ is finite. Fix an orbital digraph
$\Gamma$ of $G$, and a vertex $\alpha \in \Omega$. Let $b_r(\alpha,
\Gamma)$ be the ball-size $|B_r(\alpha, \Gamma)|$, and let
$N_r(\alpha, \Gamma)$ be the number of suborbits of $G_\alpha$ in
$B_r(\alpha, \Gamma)$. The sequence $(b_r/N_r)$ is then the {\em
average subdegree sequence} of $G$ with respect to the digraph $\Gamma$.
The growth of this sequence can be used as a measure of the subdegree
growth of $G$.

The usefulness of average subdegree growth as a measure of subdegree
growth is limited, however, by its dependence on the orbital digraph
chosen. While it is in fact possible to explicitly bound this
dependence, a further, and more serious, limitation is an inherent
lack of subtlety. As an illustration, suppose we are given integers
$m, t \geq 2$. The groups $G:= \aut \Gamma(m, K_{t+1})$ and $H:=G \Wr
\sym(2)$ acting on $\Omega := V\Gamma(m, K_{t+1})$ and $\Omega \times
\Omega$ respectively are manifestly dissimilar: all non-diagonal
orbital digraphs of the former have infinitely many ends; all
non-diagonal orbital digraphs of the latter have just one. Their many
differences are encapsulated by their lower subdegree growth rates,
with $G$ exhibiting exponential growth, and $H$ subexponential
non-polynomial growth. However, if one were instead to compare the
average subdegree growth rates the groups, no such distinction is
possible.

Indeed, define $\Gamma := \Gamma(m, K_{t+1})$ and $\Lambda :=
(\Omega^2, ((\alpha, \alpha), (\alpha, \beta))^H)$, where $\alpha$ and
$\beta$ are adjacent in $\Gamma$. The average subdegree growth rate of
$G$ with respect to $\Gamma$ is exponential, with
\begin{align*}
\lim_{r \rightarrow \infty} \left ( \frac{|B_r(\alpha,
\Gamma)|}{N_r(\alpha, \Gamma)} \right )^{1/r}
&\geq \, \left ( \frac{m(m-1)^{r-1}t^r}{r} \right )^{1/r}\\
&= \, (m-1)t.
\end{align*}
Furthermore, $|B_r(\alpha, \Gamma)| \leq r|S_r(\alpha, \Gamma)|$, so
the limit is equal to $(m-1)t$.

By Lemma~\ref{lemma:wreath_prod_graph_sphere},
\[|S_r((\alpha, \alpha), \Lambda)| = \sum_{k=0}^r |S_k(\alpha,
\Gamma)| |S_{r-k}(\alpha, \Gamma)|,\]
so
\[|B_r((\alpha, \alpha), \Lambda)| \leq \frac{m (r+1) |B_r(\alpha,
\Gamma)|}{(m-1)}.\]
Hence $|B_r(\alpha, \Gamma)| \leq |B_r((\alpha, \alpha), \Lambda)|
\leq K (r+1) |B_r(\alpha, \Gamma)|$, where $K := m/(m-1)$.

Two vertices $(\gamma_1, \gamma_2), (\delta_1, \delta_2) \in V\Lambda$
lie in the same $H_{(\alpha, \alpha)}$ orbit if and only if the sets
$\{d_\Gamma(\alpha, \gamma_1), d_\Gamma(\alpha, \gamma_2)\}$ and
$\{d_\Gamma(\alpha, \delta_1), d_\Gamma(\alpha, \delta_2)\}$ are
equal. Thus if $n_r$ denotes the number of $(\alpha,
\alpha)$-suborbits in $S_r((\alpha, \alpha), \Lambda)$ then $1 \leq
n_r \leq r$, so
\[r \leq N_r((\alpha, \alpha), \Lambda) \leq r^2.\]
Whence the average subdegree growth of $H$ is exponential, with
\[\lim_{r \rightarrow \infty} \left ( \frac{|B_r(\alpha,
\Gamma)|}{N_r(\alpha, \Gamma)} \right )^{1/r} = (m-1)t.\]

In light of these problems we will henceforth focus only on the
relationship between the lower and upper subdegree growth rates of a
group, and its structure.

\section{Subdegree growth and group structure}

In this section we give bounds on the growth of the lower and upper
subdegree sequences of infinite locally finite primitive groups, and
show that, when the growth of the lower subdegree sequences is fast
enough, the rate of growth uniquely determines the group. Following
this, we show that in many cases the number of permutation ends of a
group is determined by its multiset of subdegrees.

\begin{thm} {\normalfont(\cite{praeger})} \label{thm:praeger} Let
$\Gamma$ be an infinite connected vertex- and arc-transitive
digraph with finite but unequal in-valency and out-valency. Then there
is an epimorphism $\varphi$ from the vertex set of $\Gamma$ to the set
of integers $\mathbb{Z}$ such that $(\alpha, \beta)$ is an arc of
$\Gamma$ only if $\varphi(\beta) = \varphi(\alpha)+1$. \qed \end{thm}

\begin{cor} \label{cor:paired_suborbits_have_same_size} Suppose $G$ is
a primitive group of permutations of an infinite set $\Omega$, and
every suborbit of $G$ is finite. If $\alpha \in \Omega$ and
$\Delta(\alpha)$ and $\Delta^*(\alpha)$ are paired $\alpha$-suborbits
then $|\Delta(\alpha)| = |\Delta^*(\alpha)|$. \end{cor}

\begin{proof} Suppose $G$ is a primitive group of permutations of an
infinite set $\Omega$, and every suborbit of $G$ is finite. Fix an
element $\alpha \in \Omega$ and let $\Delta(\alpha)$ and
$\Delta^*(\alpha)$ be paired $\alpha$-suborbits. Suppose
$|\Delta(\alpha)| \not = |\Delta^*(\alpha)|$. We will show this
implies $G$ cannot be primitive, contradicting our original
assumption.

Choose $\beta \in \Delta(\alpha)$ and let $\Gamma$ be the orbital
digraph $(\Omega, (\alpha, \beta)^G)$. As $G$ is primitive, this digraph
is connected, and because $\Delta(\alpha)$ and $\Delta^*(\alpha)$ are
finite, $\Gamma$ is locally finite, with in-valency $|\Delta(\alpha)|$
and out-valency $|\Delta^*(\alpha)|$. As these are not equal, we may
apply Theorem~\ref{thm:praeger} to deduce that there is an epimorphism
$\varphi$ from the vertex set of $\Gamma$ to the set of integers
$\mathbb{Z}$ such that $(\gamma, \delta)$ is an arc of $\Gamma$ only
if $\varphi(\delta) = \varphi(\gamma)+1$.

Without loss of generality, we may assume $\varphi(\alpha) = 0$. If
$\alpha_0 \alpha_1 \cdots \alpha_n$ is any cycle in $\Gamma$ with
$\alpha_0 = \alpha$, we must have $\sum_{i=0}^n \varphi(\alpha_i) =
0$. Thus, there can be no odd cycles in $\Gamma$ containing $\alpha$.
Since $\Gamma$ is vertex-transitive, it contains no odd cycles. Hence
the $G$-congruence given by $\gamma \cong \delta$ if and only if
$d_\Gamma(\gamma, \delta)$ is even, is non-trivial and non-universal.
Whence, $G$ is not primitive. \end{proof}

Let $G$ be a primitive group of permutations of an infinite set
$\Omega$, with every subdegree of $G$ finite, and suppose $\Gamma$ is
an orbital digraph of $G$. Fix $\alpha \in \Omega$. An upper bound can
be found for the growth of the lower subdegree sequences by bounding
the growth of $s_r:=|S_r(\alpha, \Gamma)|$, since $G_\alpha$ fixes
$S_r(\alpha, \Gamma)$ setwise.

\begin{lemma} \label{lemma:growth_leq_s_1_minus_1} If $\Gamma$ is an
infinite locally finite primitive digraph then, for all $r \geq 1$,
\[ s_r \leq s_1 (s_1 - 1)^{r-1}. \] \end{lemma}

\begin{proof} Since $\Gamma$ is vertex-primitive, it is connected and
every vertex has valency $s_1$. Thus for all $r > 1$, any vertex in
the sphere $S_r(\alpha, \Gamma)$ is connected to at least one vertex
in $S_{r-1}(\alpha, \Gamma)$, and is therefore adjacent to at most
$s_1-1$ vertices in $S_{r+1}$. Whence $s_r \leq s_{r-1}
(s_1-1)$.\end{proof}

A corresponding bound is easily obtained for the growth of the upper
subdegree sequence.

\begin{lemma} Let $G$ be a locally finite primitive group of
permutations of an infinite set $\Omega$. If $(M_r)$ is the upper
subdegree sequence of $G$, then for all $r \geq 1$,
\[M_r \leq 2M_1 (2M_1-1)^{r-1}.\]\end{lemma}

\begin{proof} Fix $\alpha \in \Omega$ and choose $\beta \in \Omega$
such that the suborbit $\beta^{G_\alpha}$ is of size $M_1$. Let
$\Gamma$ be the orbital digraph $(\Omega, (\alpha, \beta)^G)$, and let
$\Delta(\alpha)$ denote the suborbit $\beta^{G_\alpha}$.

If $\Delta^*(\alpha)$ is the suborbit paired with $\Delta(\alpha)$
then $|\Delta^*(\alpha)|$ is also equal to $M_1$ by
Corollary~\ref{cor:paired_suborbits_have_same_size}. Since $S_1 =
\Delta(\alpha) \cup \Delta^*(\alpha)$, the valency of $\Gamma$ is at
most $2M_1$.

If the upper subdegree sequence contains just one subdegree then there
is nothing to prove, so suppose this is not the case. Choose $r \geq
1$ such that $M_r$ and $M_{r+1}$ are elements of the upper subdegree
sequence of $G$; since all members of this sequence are distinct,
there exists an integer $t \geq 1$ such that $M_r$ is the largest
subdegree in the sphere $S_t$ but not in $S_{t+1}$. Thus, there exists
a vertex $\gamma \in S_t$ such that $\gamma$ is adjacent to a vertex
$\delta \in S_{t+1}$ with $|\delta^{G_\alpha}| > M_r$.

Since $|\delta^{G_\alpha}| > M_r$, we have $|\delta^{G_\alpha}| \geq
M_{r+1}$, and therefore
\[M_{r+1} \leq |\delta^{G_\alpha}| \leq (s_1-1)|\gamma^{G_\alpha}|
\leq M_r(s_1-1). \qedhere \]
\end{proof}

If $\Gamma$ is a locally finite primitive digraph with connectivity one,
and $T$ is the block-cut-vertex-tree of $\Gamma$, we define the {\it
lobe-distance} between two vertices $\alpha$ and $\beta$ in $\Gamma$
to be
\[bd(\alpha, \beta) := \frac{d_T(\alpha, \beta)}{2}.\]
Informally, one may think of the lobe-distance between $\alpha$ and
$\beta$ as being the number of lobes of $\Gamma$ through which any
geodesic between $\alpha$ and $\beta$ must pass.

If $\Gamma$ is an infinite locally finite primitive digraph and $\alpha
\in V\Gamma$, we define $n_r(\Gamma)$ to be the number of
$\alpha$-suborbits of $\aut \Gamma$ in $S_r(\alpha, \Gamma)$. Since
$G$ acts transitively on $\Gamma$, this definition is independent of
our choice of $\alpha$. When there can be no ambiguity as to the
identity of the digraph in question, $n_r(\Gamma)$ will be written as
$n_r$.

Let $\{f_r\}_{r \geq 0}$ be the sequence of {\em Fibonacci numbers}
defined by the relation $f_{r+1} = f_r + f_{r-1}$ for $r \geq 1$ and
$f_0 = f_1 = 1$.

\begin{thm} If $\Gamma$ is a locally finite primitive digraph with
connectivity one that is not distance-transitive then $n_r \geq f_r$
for all $r \geq 1$.\end{thm}

\begin{proof} Fix $\alpha \in V\Gamma$. To each vertex $\beta$ in
$V\Gamma$ choose a geodesic in $\Gamma$ from $\alpha$ to $\beta$ and
assign a label $\ell(\beta):=(l_1, \ldots, l_k)$, where $k =
bd(\alpha, \beta)$, and for all $r$ satisfying $1 \leq r \leq k$ the
integer $l_r$ is the number of vertices at lobe-distance $r$ from
$\alpha$ in the geodesic. It is simple to check that, since $\Gamma$
has connectivity one, this is independent of the geodesic chosen. Let
$L_r$ be the set of labels for vertices in $S_r(\alpha, \Gamma)$ and
put $k_r := |L_r|$. It is very easy to see two vertices $\beta,
\beta^\prime \in S_r(\alpha, \Gamma)$ lie in the same
$\alpha$-suborbit only if $\ell(\beta) = \ell(\beta^\prime)$. Thus
$n_r(\Gamma) \geq k_r$.

By Theorem~\ref{thm:classification_of_dist_trans_graphs}, if $\Gamma$
is not distance-transitive then the lobes of $\Gamma$ have diameter at
least $2$. We claim that, in this case,
\begin{equation} \label{equation:no_label_in_k1_graphs} k_{r+1} \geq
k_r + k_{r-1}. \end{equation}
Indeed, for all $r \geq 1$ there is an injective correspondence from
$L_r$ into $L_{r+1}$ via the map sending $(l_1, \ldots, l_k) \in L_r$
to $(l_1, \ldots, l_k, 1) \in L_{r+1}$. Fix $r \geq 1$ and note there
are $k_{r-1}$ labels in $L_{r}$ whose last entry is $1$. If $(l_1,
\ldots, l_{k-1}, 1)$ is such a label, then $(l_1, \ldots, l_{k-1}, 2)$
is a label in $L_{r+1}$. Therefore, there are at least $k_{r-1}$
labels in $L_{r+1}$ whose final entry is $2$; we have already seen
there are $k_r$ labels in $L_{r+1}$ whose final entry is $1$, so we
must have $k_{r+1} \geq k_r + k_{r-1}$, as claimed. Since $k_0 = k_1 =
1$ we have $n_r \geq k_r \geq f_r$ for all $r \geq 0$.\end{proof}

It is well known that
\[\lim_{r \rightarrow \infty} f_r^{1/r} = \frac{1 + \sqrt{5}}{2};\]
consequently it is possible to find a lower bound for the growth of
the sequence $(n_r^{1/r})$.

\begin{cor} \label{thm:connectivity1_suborbit_growth} If $\Gamma$ is a
locally finite primitive digraph with connectivity one that is not
distance-transitive then
\[\liminf_{r \rightarrow \infty} n_r^{1/r} \geq \frac{1 + \sqrt{5}}{2}. \qedhere \] \end{cor}

The bound given in the above corollary is sharp; that is, it cannot be
improved upon. Indeed, consider the Peterson graph $P_5$. This is a
finite primitive distance-transitive graph with diameter $2$.
It can be considered to be directed by inserting arcs, one in each direction,
between adjacent vertices. Let
$\Gamma:=\Gamma(2, P_5)$ and $G:= \aut \Gamma$, and fix $\alpha \in
V\Gamma$. It is hopefully clear that two vertices $\beta, \beta^\prime
\in S_r(\alpha, \Gamma)$ lie in the same orbit of $G_\alpha$ if and
only if $\ell(\beta) = \ell(\beta^\prime)$. Thus $n_r = k_r$. We
claim, for all $r \geq 1$,
\begin{equation} \label{equation:labels_for_Gamma_P5} k_{r+1} = k_r +
k_{r-1}. \end{equation}
Indeed, let $L_r^\prime$ be the elements of $L_r$ whose last entry is
$1$, let $L_r^{\prime \prime}$ be those elements whose last entry is
$2$ and put $k_r^\prime:=|L_r^\prime|$ and $k_r^{\prime
\prime}:=|L_r^{\prime \prime}|$. Since $P_5$ has diameter $2$ we have
$k_r = k_r^\prime + k_r^{\prime \prime}$. We may assign to each label
$\ell \in L_r^{\prime \prime}$ the unique label in $L_{r-1}^\prime$
obtained by changing the last entry in $\ell$ from $2$ to $1$. Hence
$k_r'' \geq k_{r-1}'$. We may also assign to each label $\ell \in
L_{r-1}^\prime$ a unique label in $L_r^{\prime \prime}$ obtained by
changing the last entry in $\ell$ to $2$, so $k_{r-1}' \geq k_r''$.
Whence $k_r^{\prime \prime} = k_{r-1}^\prime$. Similarly, we assign to
each label $(l_1, \ldots, l_k) \in L_{r-1}$ the label $(l_1, \ldots,
l_k, 1) \in L_r^\prime$, and to each label $(l_1, \ldots, l_k, 1) \in
L_r^\prime$ the label $(l_1, \ldots, l_k) \in L_{r-1} $. Such a
correspondence is bijective, so $k_{r-1} = k_r^\prime$. Thus $k_{r+1}
= k_{r+1}^\prime + k_{r+1}^{\prime \prime} = k_{r+1}^\prime +
k_r^\prime = k_r + k_{r-1}$, as claimed.

Now $k_0 = k_1 = 1$ and so by (\ref{equation:labels_for_Gamma_P5}),
for all $r \geq 0$ we have $n_r = k_r = f_r$. Thus
\[ \lim_{r \rightarrow \infty} n_r^{1/r} = \frac{1 + \sqrt{5}}{2}.\]

Let $\Gamma$ be a infinite locally finite connected vertex- and
arc-transitive digraph, and fix $\alpha \in V\Gamma$. For $\gamma \in
S_r(\alpha)$ we define
\begin{align*}
a(\gamma) :=& |S_1(\alpha) \cap S_{r}(\gamma)|;\\
b(\gamma) :=& |S_1(\alpha) \cap S_{r+1}(\gamma)|;\\
c(\gamma) :=& |S_1(\alpha) \cap S_{r-1}(\gamma)|. \end{align*}
The following lemma is an extension of an observation by Macpherson in
\cite{macpherson:dtg}. The argument presented here is based on that
given by Dicks and Dunwoody in \cite{dicks&dunwoody:gaog}.

\begin{lemma} \label{lemma:constant_c_means_more_then_one_end} If
there exists a natural number $R_0$ such that, for all $\beta, \gamma
\in V\Gamma$ with $d(\alpha, \beta) > R_0$ and $d(\alpha, \gamma) >
R_0$, we have $c(\beta) = c(\gamma)$ and $b(\beta) = b(\gamma)$, then
$\Gamma$ has more than one end. \end{lemma}

\begin{proof}
Suppose, for all $\beta, \gamma \in V\Gamma$ with $d(\alpha, \beta) >
R_0$ and $d(\alpha, \gamma) > R_0$, we have $c(\beta) = c(\gamma)$ and
$b(\beta) = b(\gamma)$. We will describe an infinite set of vertices
$s$, such that the complement $s^*:=V\Gamma \setminus s$ is infinite,
and the set of arcs $\delta s$ from $s$ to $s^*$ is finite, thus
showing $\Gamma$ must have more than one end.

Fix an arc $(\alpha_0, \alpha_1) \in V\Gamma$ and define $s:=\{\gamma
\in V\Gamma \mid R_0+1 \leq d(\alpha_0, \gamma) = d(\alpha_1,
\gamma)+1\}$. This set and its compliment are both infinite. Indeed,
given a positive integer $n$ one may choose a vertex $\gamma \in
V\Gamma$ with $d(\alpha_0, \gamma) = 2n+1$. There exists a geodesic
$\alpha_0 \beta_1 \beta_2 \ldots \beta_{2n} \gamma$ of length $2n+1$
between $\alpha_0$ and $\gamma$. If $e$ is the arc between $\beta_n$
and $\beta_{n+1}$ then, since $\Gamma$ is arc-transitive, there
exists an automorphism of $\Gamma$ mapping $e$ to the arc $(\alpha_0,
\alpha_1)$. Therefore, both $s$ and $s^*$ contain vertices in
$S_n(\alpha_0, \Gamma) \cup S_{n+1}(\alpha_0, \Gamma)$. Whence $s$ and
$s^*$ are infinite.

We now show $\delta s$ is finite. Suppose $e \in E$ is an arc between
$\beta \in s$ and $\gamma \in s^*$. Write $i = d(\alpha_0, \beta)$.

We claim $d(\alpha_0, \gamma) = 1+d(\alpha_1, \gamma)$. Observe
$d(\alpha_0, \gamma) \leq d(\alpha_0, \beta) + d(\beta, \gamma) = i +
1$, and $d(\alpha_0, \beta) \leq d(\alpha_0, \gamma) + d(\beta,
\gamma)$, so $d(\alpha_0, \gamma) \geq d(\alpha_0, \beta) - d(\beta,
\gamma) = i - 1$. Hence,
\[i-1 \leq d(\alpha_0, \gamma) \leq i+1.\]
We now consider three cases: when $d(\alpha_0, \gamma) = i-1$, when
$d(\alpha_1, \gamma) = i$ and finally when $d(\alpha_0, \gamma) \geq
i$ and $d(\alpha_1, \gamma) \leq i-1$, and in each case show the claim
is true.

Suppose $d(\alpha_0, \gamma) = i-1$. Then $d(\alpha_0, \beta) =
d(\alpha_0, \gamma) + d(\gamma, \beta)$, so there is a geodesic from
$\alpha_0$ to $\beta$ that contains $\gamma$. Hence
\[ S_1(\alpha_0) \cap S_{i-2}(\gamma) \subseteq S_1(\alpha_0) \cap
S_{i-1}(\beta);\]
however, since $i-1 \geq R_0$ we have $|S_1(\alpha_0) \cap
S_{i-2}(\gamma)| = c(\gamma) = c(\beta) = |S_1(\alpha_0) \cap
S_{i-1}(\beta)|$, so the two sets must be equal. Now $\alpha_1 \in
S_1(\alpha_0) \cap S_{i-1}(\beta)$, so $\alpha_1 \in S_1(\alpha_0)
\cap S_{i-2}(\gamma)$; that is, $d(\alpha_1, \gamma) = i-2 =
d(\alpha_0, \gamma) - 1$ as claimed.

Next, suppose $d(\alpha_1, \gamma) = i$. Then $d(\alpha_1, \gamma) =
d(\alpha_1, \beta) + d(\gamma, \beta)$, so there is a geodesic from
$\alpha_1$ to $\gamma$ that contains $\beta$. Hence
\[ S_1(\alpha_1) \cap S_{i+1}(\gamma) \subseteq S_1(\alpha_1) \cap
S_{i}(\beta);\]
however, since $i-1 \geq R_0$ we have $|S_1(\alpha_1) \cap
S_{i+1}(\gamma)| = b(\gamma) = b(\beta) = |S_1(\alpha_1) \cap
S_{i}(\beta)|$, so the two sets must be equal. Now $\alpha_0 \in
S_1(\alpha_1) \cap S_{i}(\beta)$, so $\alpha_0 \in S_1(\alpha_1) \cap
S_{i+1}(\gamma)$; that is, $d(\alpha_0, \gamma) = i+1 = d(\alpha_1,
\gamma) + 1$ as claimed.

Finally, if $d(\alpha_0, \gamma) \geq i$ and $d(\alpha_1, \gamma) \leq
i-1$, then $d(\alpha_1, \gamma) + 1 \leq i \leq d(\alpha_0, \gamma)$.
In fact, since $d(\alpha_0, \gamma) \leq d(\alpha_1, \gamma) +
d(\alpha_1, \alpha_0) = d(\alpha_1, \gamma) + 1$, we have $d(\alpha_0,
\gamma) \geq 1+d(\alpha_1, \gamma)$, so $d(\alpha_0, \gamma) = 1 +
d(\alpha_1, \gamma)$ as claimed.

Now $\gamma \in s^*$, so either $d(\alpha_0, \gamma) < R_0+1$ or
$d(\alpha_0, \gamma) \not = d(\alpha_1, \beta)+1$. The latter is not
true by the above argument, so we must have $d(\alpha_0, \gamma) <
R_0+1$. Hence, there are only finitely many $\gamma \in s^*$ that are
adjacent to a vertex in $s$, so $\delta s$ is finite, and $\Gamma$ has
more than one end. \end{proof}

\begin{thm} \label{thm:one_end_suborbit_growth} If $\Gamma$ is an
infinite vertex- and arc-transitive locally finite digraph with one end
then $n_r(\Gamma) \geq 2$ for all large enough $r$. \end{thm}
\begin{proof} Let $\Gamma$ be an infinite locally finite digraph with
one end, and fix $\alpha \in V\Gamma$. We claim there is an integer
$R_0$ such that for all $r \geq R_0$ there exist vertices $\gamma_1,
\gamma_2 \in S_r(\alpha)$ with $c(\gamma_1) \not = c(\gamma_2)$ or
$b(\gamma_1) \not = b(\gamma_2)$.

Suppose no such $R_0$ exists. Then there exists an infinite sequence
$(r_i)$ such that for all $\gamma_1, \gamma_2 \in S_{r_i}(\alpha)$
we have $c(\gamma_1) = c(\gamma_2)$ and $b(\gamma_1) = b(\gamma_2)$.
For each $i \geq 1$ choose $\gamma \in S_{r_i}(\alpha)$ and set
$c_{r_i}:=c(\gamma)$ and $b_{r_i}:=b(\gamma)$. It is easy to see
$c_{r_i} \geq c_{r_{i-1}}$ and $b_{r_i} \leq b_{r_{i-1}}$. Since $1
\leq c_{r_i} \leq |S_1(\alpha)|$ and $1 \leq b_{r_i} \leq
|S_1(\alpha)|$ for all $i \geq 0$, there exists constants $k, c$ and
$b$ such that, for all $i \geq k$ we have $c_{r_i} = c$ and $b_{r_i} =
b$. Put $R_0 := r_k$. Suppose there exists $r \geq R_0$ such that
$S_r(\alpha)$ contains two vertices $\gamma_1$ and $\gamma_2$ with
$c(\gamma_1) \not = c(\gamma_2)$. Without loss of generality, one may
suppose $c(\gamma_1) > c(\gamma_2)$. However, $r > R_0$ so
$c(\gamma_2) \geq c$ and $c(\gamma_1) > c$. If $i$ is chosen so $r_i >
r$ then $c_{r_i} \geq c(\gamma_1) > c$ which is a contradiction.

It must therefore be the case that, for all $r \geq R_0$, we have
$c(\gamma) = c$ for all $\gamma \in S_r(\alpha)$. A similar argument
shows that for all $r \geq R_0$, we have $b(\gamma) = b$ for all
$\gamma \in S_r(\alpha)$. Hence, by
Lemma~\ref{lemma:constant_c_means_more_then_one_end}, the digraph
$\Gamma$ must have more than one end. Since this is not the case, our
claim must be true.

Let $G:=\aut \Gamma$. It is clear that two vertices $\gamma_1,
\gamma_2 \in S_r(\alpha)$ lie in the same orbit of $G_\alpha$ only if
$c(\gamma_1) = c(\gamma_2)$ and $b(\gamma_1) = b(\gamma_2)$. Hence,
for all $r > R_0$ we have $n_r(\Gamma) \geq 2$.\end{proof}

\begin{thm} \label{thm:subdegree_growth} Let $G$ be a primitive group
of permutations of an infinite set $\Omega$. If $G$ is locally finite
with more than one permutation-end then the lower subdegree sequence
of $G$ grows exponentially if and only if $G$ is distance-transitive.
In this case, $G$ has height $\omega$, and its subdegree sequence
$(m_r)$ satisfies
\[\lim_{r \rightarrow \infty} m_r^{1/r} = \left ( \frac{m-1}{m} \right ) m_1,\]
where $G$ acts distance-transitively on the distance-transitive
infinite locally finite digraph $\Gamma(m, K_{t+1})$. Furthermore, if the
growth of the lower subdegree sequence of $G$ is not exponential, then
it is polynomial. \end{thm}

\begin{proof} Suppose $G$ is a primitive group of permutations of an
infinite set $\Omega$, possessing an orbital digraph with more than one
end, and $G$ has a finite suborbit whose pair is also finite. By
Theorem~\ref{thm:finite_suborbit_implies_all_suborbits_finite}, every
suborbit of $G$ is finite. Let $(m_r)$ be the lower subdegree sequence
of $G$. Note that, if $G$ acts distance-transitively on the locally
finite infinite distance-transitive digraph $\Gamma(m, K_{t_1 +1})$,
then every suborbit is self-paired, $G$ has height $\omega$, so
$(m_r)$ is equal to the subdegree sequence of $G$, and the subdegree
growth of $G$ is exponential with
\begin{align*} \lim_{r \rightarrow \infty} m_r^{1/r}
    =& \left ( \frac{m-1}{m} \right ) m t_1\\
    =& \left ( \frac{m-1}{m} \right ) m_1. \end{align*}

Now consider the converse. Suppose the group $G$ does not act
distance-transitively on any locally finite infinite
distance-transitive digraph. We will show that the lower subdegree
growth of $G$ is bounded above by some polynomial.

By Corollary~\ref{thm:connectivity1_suborbit_growth}, $\liminf
n_r^{1/r} \geq (1+\sqrt{5})/2 > 3/2$. Hence, there exists an integer
$R$ such that, for all $r > R$, we have $n_r > (3/2)^r$, and thus $N_r
> (3/2)^r$.

Fix $r > R$ and observe that $m_{N_r} \leq s_1(s_1-1)^{r-1}$ by
Lemma~\ref{lemma:growth_leq_s_1_minus_1}. We may choose an integer $N$
such that, for all $n \geq N$, we have $(3/2)^{n} \geq s_1(s_1-1)$.
Thus, for all $n \geq N$,
\[m_{N_r} < N_r^n.\]
Furthermore, given any integer $s$ with $N_{r-1} \leq s \leq N_r$, the
subdegree $m_s$ satisfies $m_s \leq s_1(s_1-1)^{r-1}$. Since $s \geq
N_{r-1}$ we also have
$s^n \geq N_{r-1}^n > (3/2)^{n(r-1)} > s_1(s_1-1)^{r-1} \geq m_s$, so
the growth of the lower subdegree sequence of $G$ is polynomial.
\end{proof}

If one removes the condition that $G$ have more than one
permutation-end then the following is obtained.

\begin{thm} \label{thm:general_subdegree_growth} Suppose $G$ is an
infinite locally finite group of permutations of an infinite set
$\Omega$, and $G$ is not distance-transitive. If $(m_r)$ is the lower
subdegree sequence of $G$ then
\[\liminf_{r \rightarrow \infty} m_r^{1/r} \leq \sqrt{2 m_1 - 1}.\]\end{thm}

\begin{proof} Suppose $G$ is a locally finite primitive group of
permutations of an infinite set $\Omega$ and does not act
distance-transitively on any digraph. If $G$ has an orbital digraph with
more than one end then its lower subdegree growth is subexponential by
Theorem~\ref{thm:subdegree_growth}. So, suppose $G$ has an orbital
digraph with precisely one end; every orbital digraph of $G$ therefore has
one permutation-end. Let $\Delta(\alpha)$ be a suborbit of size $m_1$,
and let $\Gamma$ be the orbital digraph $(\Omega, \Delta)$. Then
$S_1(\alpha, \Gamma) = \Delta(\alpha) \cup \Delta^*(\alpha)$. Let
$s_r:=|S_r(\alpha, \Gamma)|$ and $n_r := n_r(\Gamma)$, and let $N_r$
be the sum $\sum_{i=1}^r n_i$. Since all suborbits of $G$ are finite,
one may deduce from
Corollary~\ref{cor:paired_suborbits_have_same_size} that
$\Delta^*(\alpha)$, the suborbit paired with $\Delta(\alpha)$, also
has cardinality $m_1$. Thus, $s_1 \leq 2 m_1$.

From the proof of Theorem~\ref{thm:subdegree_growth}, $m_{N_r} \leq
s_1 (s_1-1)^{r-1}$. We again note it is sufficient to show
this result holds when
$G=\aut \Gamma$.

By Theorem~\ref{thm:one_end_suborbit_growth}, there exists an integer
$R_0$ such that $n_r \geq 2$ for all $r \geq R_0$. Hence, for
sufficiently large $r$,
\[N_r \geq R_0 + 2(r-R_0).\]
Thus, if $a > 1$,
\begin{align*} \limsup_{r \rightarrow \infty} m_{N_r}^{1/N_r}
    \leq& \limsup_{r \rightarrow \infty} \, (s_1 - 1)^{r/N_r}\\
    \leq& \limsup_{r \rightarrow \infty} \, (s_1 - 1)^{r/(2r - R_0)}\\
    =& \sqrt{2 m_1-1}. \end{align*}
Hence $\displaystyle{\liminf_{r \rightarrow \infty} m_r^{1/r} \leq
\sqrt{2 m_1-1}}$ as required. \end{proof}

\begin{cor} If $G$ is a group acting primitively on an infinite set
$\Omega$ with a finite suborbit whose pair is also finite, then the
subdegrees of $G$ are all finite. If $(m_r)$ is the lower subdegree
sequence of $G$ then
\[\liminf_{r \rightarrow \infty} m_r^{1/r} > \sqrt{2 m_1 - 1}\]
if and only if $G$ acts distance-transitively on some distance-transitive locally finite
infinite digraph $\Gamma(m, K_{t+1})$ with $m > 2$
and $t \geq 2$, or $m=2$ and $t \geq 4$. \end{cor}

\begin{proof} Suppose $G$ is a locally finite primitive group of
permutations of an infinite set $\Omega$. If $G$ does not act
distance-transitively on any locally finite orbital digraph and $(m_r)$
is the lower subdegree sequence of $G$ then
\[\liminf_{r \rightarrow \infty} m_r^{1/r} \leq \sqrt{2 m_1 - 1}\]
by Theorem~\ref{thm:general_subdegree_growth}.

Now suppose $G$ acts distance-transitively on a locally finite
distance-transitive digraph $\Gamma$. By
Theorem~\ref{thm:classification_of_dist_trans_graphs}, we may write
$\Gamma = \Gamma(m, K_{t+1})$ for some $m \geq 2$ and $t \geq 2$.
Observe that ${\displaystyle\liminf_{r \rightarrow \infty} m_r^{1/r} =
(m-1)t}$, so the limit ${\displaystyle \liminf_{r \rightarrow \infty}
m_r^{1/r} > \sqrt{2 m_1 - 1}}$ precisely when $m > 2$ and $t \geq 2$,
or $m=2$ and $t \geq 4$. \end{proof}

All known examples of primitive groups with locally finite one-ended
orbital digraphs exhibit subexponential lower subdegree growth;
furthermore, it seems highly unlikely that examples exhibiting
exponential growth exist.

\begin{conjecture} \label{conjecture:final_conjecture} If $G$ is a
group acting primitively on an infinite set $\Omega$ with a finite
suborbit whose pair is also finite, then the subdegrees of $G$ are all
finite and $G$ has exponential lower subdegree growth if and only if
$G$ is distance-transitive. \end{conjecture}

Of course there are many further questions that remain unanswered. Which primitive groups
exhibit subexponential non-polynomial subdegree growth, and are there
gaps in growth rates that allow one to determine a group given
its subdegree growth rate? What is the relationship between the
permutation-ends of a group and its subdegree growth rates? This latter
question is the focus of the final section of this paper. In it, we detail
some preliminary results that illustrate a relationship does indeed exist; however,
much work remains to be done before its nature is fully determined.

\section{Subdegree growth and ends of orbital digraphs}

We begin by observing that the relationship between the subdegree
growth of a primitive group $G$ and its permutation-end structure is
more subtle than one might expect.

\begin{thm} If $(m_r)$ is the lower subdegree sequence of a locally
finite infinite primitive group $G \leq \sym(\Omega)$, then there exist
infinite primitive groups $G'$ and $G''$ whose suborbits are all
finite, with lower subdegree sequences $(m_r')$ and $(m_r'')$
respectively, such that $G'$ has one permutation-end, and $G''$ has
infinitely many permutation-ends, with
\[m_r' \leq 2 m_r\]
and
\[m_r'' \leq 2 m_r\]
for all $r \geq 1$. \end{thm}

\begin{proof} Fix $\alpha \in \Omega$, and write
$\underline{\alpha}:=(\alpha, \alpha) \in \Omega^2$. Take $G'$
to be the wreath product $G \Wr \sym(2)$, and consider its
product action on $\Omega^2$. By
Theorem~\ref{thm:wreath_product_construction}, all suborbits of $G'$
are finite and every orbital digraph has one end. Furthermore, for every
suborbit $\beta^{G_\alpha}$ of $G$, the set $(\beta,
\alpha)^{G_{\underline{\alpha}}'}$ is a suborbit of $G'$; since
$|(\beta, \alpha)^{G'_{\underline{\alpha}}}| = 2 |\beta^{G_\alpha}|$,
we have $m_r' \leq 2 m_r$.

Let $G''$ be the group $G(2, G)$ constructed in
Chapter~\ref{section:constructions_with_infinitely_many_ends}, and let
$\Lambda$ be an orbital digraph of $G$ acting on $\Omega$ with
connectivity greater than one; by
Theorem~\ref{thm:infinitely_ended_construction}, $G''$ acts
primitively on the vertex set of the digraph $\Gamma(2, \Lambda)$.
Indeed, this digraph is an orbital digraph of $G''$. Since $\Gamma(2,
\Lambda)$ has infinitely many ends, every orbital digraph of $G''$ has
infinitely many permutation-ends. Furthermore, the action of
$G''_{\{\Lambda\}}$ on $V\Lambda$ is isomorphic to the action of $G$
on $V\Lambda = \Omega$; thus, for each subdegree $m_r$ of $G$, the
group $G_{\underline{\alpha}, \{\Lambda\}}''$ has an orbit on
$V\Lambda$ of size $m_r$. Hence, for each subdegree $m_r$ of $G$, the
group $G''$ has a suborbit of size $2 m_r$. Whence, $m_r'' \leq 2
m_r$. \end{proof}

It should be noted that both the difference between $m_r$ and $m_r'$,
and the difference between $m_r$ and $m_r''$, may grow arbitrarily
large, leaving gaps in the range of possible rates of growth.

Certain growth rate are only exhibited by groups with precisely one
permutation-end. The following is immediate from
Theorem~\ref{thm:subdegree_growth}.

\begin{thm} If $G$ is a locally finite primitive group of permutations
of an infinite set $\Omega$, and the lower subdegree growth of $G$ is
subexponential but not polynomial, then $G$ has precisely one
permutation-end. \qed
\end{thm}

Example~\ref{ex:nonPolyNonExp} shows primitive groups with
non-polynomial but subexponential lower subdegree growth exist.
All known examples of locally finite primitive groups with exponential
lower subdegree growth have infinitely many permutation-ends. Indeed,
it seems highly likely that such growth rates cannot be achieved by
groups with just one permutation-end. If true then the above theorem
and the following conjecture would allow one to naturally partition
the non-polynomial rates of growth of infinite primitive groups
according to the number of permutation-ends possessed by each group.

\begin{conjecture} \label{conjecture:exp_iff_inf_ended} If $G$ is an infinite primitive permutation group
whose subdegrees are all finite, and the lower subdegree sequence of
$G$ grows exponentially, then $G$ has $2^{\aleph_0}$ permutation-ends.
\end{conjecture}

In Section~\ref{section:upper_and_lower_growth} it was shown that
there exist examples of infinite primitive groups with precisely one
permutation-end, and infinite primitive groups with infinitely many
permutation-ends, both possessing bounded lower subdegree sequences.
If, instead of examining just the lower subdegree sequence, one
considers the whole subdegree sequence, it is again possible to
determine the permutation-end structure of those infinite primitive
groups exhibiting specific rates of subdegree growth.\

We begin with two theorems describing the structure of primitive
groups with more than one permutation-end.

\begin{thm} {\normalfont \cite[Theorem 2.5]{me:InfPrimDigraphs}}
\label{thm:Gya_eq_Gyb_not_prim} Let $G$ be a vertex-transitive group
of automorphisms of a connectivity-one digraph $\Gamma$ whose
lobes have at least three vertices, and let $T$ be the
block-cut-vertex tree of $\Gamma$. If there exist distinct
vertices $\alpha, \beta \in V\Gamma$ such that, for some vertex
$x \in (\alpha, \beta)_T$,
\[ G_{\alpha, x} = G_{\beta, x}, \]
then $G$ does not act primitively on $V\Gamma$.\end{thm}

\begin{thm} {\normalfont \cite[Theorem 3.11]{me:OrbGraphs}}
\label{thm:connectivity1}
If $G$ is a primitive group of permutations of an infinite set
$\Omega$ with more than one permutation-end and no infinite
subdegree, then $G$ has a locally finite orbital digraph $\Gamma$
with connectivity one, whose lobes  are primitive but not
automorphism-regular, are pairwise isomorphic, have at least three
vertices and at most one end. Furthermore, if $\Lambda$ is a lobe
of $\Gamma$, then $G_{\{\Lambda\}}$ acts primitively but not
regularly on $V\Lambda$. \qed \end{thm}

Using these results, it is possible to determine precisely the
permutation-end structure of primitive groups whose subdegrees are
bounded above.

\begin{thm} If $G$ is an infinite primitive permutation group whose
subdegrees are all finite and bounded above, then $G$ has precisely
one permutation-end.
\end{thm}

\begin{proof} Suppose $G$ is an infinite primitive permutation group
with more than one permutation-end whose subdegrees are all finite. By
Theorem~\ref{thm:connectivity1}, $G$ has an orbital digraph of the form
$\Gamma(m, \Lambda)$ for some integer $m \geq 2$, and for some
primitive digraph $\Lambda$. Let $\Gamma$ denote the digraph $\Gamma(m,
\Lambda)$ and let $T$ be the block-cut-vertex tree of $\Gamma$.

Observe that if $\alpha$ and $\beta$ are vertices in $\Gamma$ and $x
\in VT$ lies on the $T$-geodesic $[\alpha, \beta]_T$ between $\alpha$
and $\beta$, then $G_{\alpha, \beta} \leq G_x$ and therefore
\[|G_{\alpha} : G_{\alpha, \beta}| = |G_{\alpha}:G_{\alpha, x}|
    |G_{\alpha, x}:G_{\alpha, \beta}|.\]
Hence, the cardinality of the suborbit $\beta^{G_\alpha}$ is equal to
the product $|x^{G_\alpha}| |\beta^{G_{\alpha, x}}|$.

Thus, if all subdegrees of $G$ are bounded above, then there exists a
finite number $k$ such that any automorphism in $G$ fixing the sphere
$S_k(\alpha, \Gamma)$ pointwise must also fix every vertex in
$\Gamma$, and therefore every vertex in $T$. However, groups with this
property cannot be primitive by Theorem \ref{thm:Gya_eq_Gyb_not_prim}.
\end{proof}

Using a similar argument, it is sometimes possible to determine the
permutation-end structure of a primitive group by knowing just one
subdegree.

\begin{thm} If $G$ is an infinite primitive permutation group whose
subdegrees are all finite, and at least one subdegree is prime, then
every orbital digraph of $G$ has precisely one end.
\end{thm}

\begin{proof} Again suppose that $G$ is an infinite primitive
permutation group with more than one permutation-end whose subdegrees
are all finite. Let $\Gamma$  be a connectivity-one orbital digraph of
$G$ of the form $\Gamma(m, \Lambda)$, the existence of which is
assured by Theorem~\ref{thm:connectivity1}, and let $T$ be the
block-cut-vertex tree of this digraph.

Recall that, given $\alpha, \beta \in V\Gamma$ and a vertex $x \in VT$
lying on the $T$-geodesic $[\alpha, \beta]_T$ between $\alpha$ and
$\beta$, the cardinality of the suborbit $x^{G_\alpha}$ is equal to
$|x^{G_\alpha}| |\beta^{G_{\alpha, x}}|$.

Since $G$ acts arc-transitively on $\Gamma$, it permutes the lobes of
$\Gamma$. By Theorem~\ref{thm:connectivity1}, the setwise stabiliser
in $G$ of each lobe acts primitively on the vertices of the lobe;
whence, for each vertex $\alpha \in V\Gamma$, the stabiliser
$G_\alpha$ transitively permutes the lobes of $\Gamma$ that contain
$\alpha$. Thus $G_\alpha$ acts transitively on the sphere $S_1(\alpha,
T)$, which has cardinality $m \geq 2$.

If $\Lambda$ is a lobe of $\Gamma$, then by
Theorem~\ref{thm:connectivity1}, $\Lambda$ has at least three
vertices, and $G_{\{\Lambda\}}$ acts primitively but not regularly on
$V\Lambda$. Thus if $\alpha \in V\Lambda$, then $G_{\alpha,
\{\Lambda\}}$ fixes no vertex in $\Lambda \setminus \alpha$. The lobe
$\Lambda$ corresponds to some vertex $x \in S_1(\alpha, T)$, so
$G_{\alpha, x}$ fixes no vertex in $S_1(x, T) \setminus \{\alpha\}$.
Therefore, for each vertex $\gamma$ in $S_1(x, T) \setminus
\{\alpha\}$, there exists a prime number $p$ dividing the cardinality
of the orbit $\gamma^{G_{\alpha, x}}$. Hence $mp$ divides
$\gamma^{G_\alpha}$.

Since $G_\alpha$ acts transitively on the sphere $S_1(\alpha, T)$, for
all vertices $\gamma$ in $S_2(\alpha, T)$ there exist primes $p$ and
$q$ such that $pq$ divides the subdegree $|\gamma^{G_\alpha}|$;
therefore, the same is true of all vertices in $T$ lying at distance
greater than two from $\alpha$. Since this includes all vertices in
$\Gamma \setminus \{\alpha\}$, no subdegree of $G$ is
prime.\end{proof}

This brief look at the
relationship between the subdegree growth rates and the
permutation-end structure of infinite locally finite primitive groups
contained an obvious omission: the class of groups that exhibit polynomial
subdegree growth. Groups with
just one permutation end, and groups with infinitely many
permutation-ends, are well represented in this class.
It would be extremely interesting to know if it is
possible to determine the
permutation-end structure of such groups from the order of their
growth.

It appears that more must be known
about the structure of locally finite graphs with precisely one end
before any significant progress in this area can be made. A clear insight into
their nature would also aid the construction of a proof of, or counterexample to, Conjecture~\ref{conjecture:exp_iff_inf_ended}. \\

Many of the results in this paper are taken from the author's DPhil thesis, completed
at the University of Oxford, under the supervision of Peter Neumann. The author would
like to thank Dr Neumann for his enthusiasm and insightful suggestions. The author
would also like to thank the EPSRC for their generous funding.


\end{document}